\documentclass{amsart}
\usepackage{amssymb}
\usepackage{amsmath}
\usepackage{amsfonts}

\theoremstyle{plain}

\numberwithin{equation}{section}
\input{tcilatex}

\begin{document}
\title{Minimax Theorems for Set-valued Maps Without Continuity Assumptions }
\author{Monica Patriche}
\maketitle

\textbf{Abstract.}{\small \ We introduce several classes of set-valued maps
with generalized}

{\small convexity }and{\small \ we obtain minimax theorems for set-valued
maps which satisfy the }

{\small introduced properties and which are not continuous. Our method
consists of the use }

{\small of a \ fixed-point theorem for weakly naturally quasi-concave
set-valued maps defined }

{\small on a simplex in a topological vector space or of a constant
selection of quasi-convex }

{\small set-valued maps.\medskip }


{\small concave set-valued map, }$S${\small -transfer }$\mu ${\small -convex
set-valued map, transfer properly }$S-${\small quasi-convex}, {\small weakly 
}$z-${\small convex set-valued map.\medskip }

\bigskip

\textbf{1. Introduction}

The classical Ky Fan inequalities [4], [5], [6] are an undeniably important
tool in the study of many important results concerning the variational
inequalities, game theory, mathematical economics, control theory and
fixed-point theory. e.g., see [1], [2], [7], [8], [10], [13]-[17],
[19]-[21], [23], [25]-[32] and the references therein. Within recent years,
many generalizations have been successfully obtained and here we must
emphasize Ky Fan's study of minimax theorems for vector-valued mappings and
for set-valued maps. We refer the reader, for instance, to Li and Wang [15],
Luo [19], Zhang and Li [31], [32], Zhang, Cheng and Li [30]. In [20], Nessah
and Tian search the condition concerning the existence of solution of
minimax inequalities for real-valued mappings, without convexity and
compactness assumptions. They define the local dominatedness property and
prove that it is necessary and further, under some mild continuity
condition, sufficient for the existence of equilibrium in minimax
inequalities. This type of characterization of the solution for minimax
theorems leads us to the question whether similar results can be obtained,
but, by keeping the convexity assumptions and by giving up the continuity
ones over the set-valued maps.

We are introduced into the extremely limited literature concerning the
minimax theorems for set-valued maps with the opportunity to see the things
from a new perspective and to propose coherent answers to the problem of the
solution existence. Our results could be particularly designed to identify
new methods of proof for this kind of problems and to assess whether the
convexity framework can be adapted to set-valued maps with two variables and
whether classes of weakened convexity can be implemented, particularly by
relying on a mechanism which takes into accont the behaviour of the maps in
the points where their values contain or not maximal (resp. minimal)
elements of certain sets of type $\tbigcup\limits_{y\in X}F(x,y)$ or $%
\tbigcup\limits_{x\in X}F(x,y).$

In this paper, we study vector minimax inequalities for set-valued maps. We
give up the condition of continuity of the set-valued maps and, instead, we
work with some new classes of generalized convexity which we introduce: $S$%
-transfer $\mu $-convexity, transfer properly $S-$quasi-convexity and weakly 
$z-$convexity. In order to prove our results, we construct a constant
selection for a quasi-convex correspondence and we use the fixed point
theorem for weakly naturally quasi-concave set-valued maps defined on a
simplex in a topological vector space (see [22]).

The article is organized as follows. In Section 2, we introduce notations
and preliminary results. In Section 3, the convex-type properties for
set-valued maps are defined and some exemples are given as well. In Section
4, we obtain two types of Ky Fan minimax inequalities for set-valued maps.
We also provide some examples to illustrate our results. Concluding remarks
are presented in Section 5.\medskip\ 

\textbf{2. Preliminaries and Notation}\smallskip \medskip

We shall use the following notations and definitions:

Let $A$ be a subset of a topological space $X.$ $2^{A}$ denotes the family
of all subsets of $A$ and $\overline{A}$ denotes the closure of $A$ in $X$.
If $A$ is a subset of a vector space, co$A$ denotes the convex hull of $A$.
If $F$, $G:X\rightrightarrows Z$ are set-valued maps, then co $G$, $G\cap
F:X\rightrightarrows Z$ are set-valued maps defined by $($co $G)(x):=$co $%
G(x)$ and $(G\cap F)(x):=G(x)\cap F(x)$ for each $x\in X$, respectively.

In this paper, we will consider $E$ and $Z$ to be real Hausdorff topological
vector spaces and we will assume that $S$ is a pointed closed convex cone in 
$Z$ with its interior int$S\neq \emptyset .$

\textbf{Definition 2.1} (see [11]). Let $A\subset Z$ be a non-empty subset.%
\newline
\ \ \ (i) A point $z\in A$ is said to be \textit{a minimal point of }$A$ iff 
$A\cap (z-S)=\{z\},$ and Min$A$ denotes the set of all minimal points of $A.$%
\newline
\ \ \ (ii) A point $z\in A$ is said to be \textit{a weakly minimal point of} 
$A$ iff $A\cap (z-$int$S)=\emptyset ,$ and Min$_{w}A$ denotes the set of all
weakly minimal points of $A.$\newline
\ \ \ (iii) A point $z\in A$ is said to be \textit{a maximal point of }$A$
iff $A\cap (z+S)=\{z\},$ and Max$A$ denotes the set of all maximal points of 
$A.$\newline
\ \ \ (iv) A point $z\in A$ is said to be \textit{a weakly maximal point of} 
$A$ iff $A\cap (z+$int$S)=\emptyset ,$ and Max$_{w}A$ denotes the set of all
weakly maximal points of $A.$

It is easy to check that Min$A\subset $Min$_{w}A$ and Max$A\subset $Max$%
_{w}A.$

\textbf{Lemma 2.1} (see [7])Let $A\subset Z$ be a non-empty compact subset.
Then, (i) Min$A\neq \emptyset ;$ (ii) $A\subset $Min$A+S;$ (iii) $A\subset $%
Min$_{w}A+$int$S\cup \{0_{F}\};$ (iv)Max$A\neq \emptyset ;$ (v) $A\subset $%
Max$A-S;$ (vi) $A\subset $Max$_{w}A-$int$S\cup \{0_{F}\}.$

\textit{Notation. If }$X$ and $Y$ are sets and $F:X\times X\rightrightarrows
Y$ is a set-valued map$,$ we will denote $F(x,X)=\tbigcup\limits_{y\in
X}F(x,y)$ and $F(X,y)=\tbigcup\limits_{x\in X}F(x,y).$

We present the following types of generalized convex mappings and set-valued
maps.

\textbf{Definition 2.2} Let $X$ be a non-empty convex subset of a
topological vector space\textit{\ }$E,$ $Z$ a real topological vector space
and $S$ a pointed closed convex cone in $Z$ with its interior int$S\neq
\emptyset .$ Let $F:X\rightrightarrows Z$ be a set-valued map with non-empty
values.

(i) $F$ is said to be (in the sense of [12 , Definition 3.6]) \textit{%
type-(iii)} \textit{properly }$S-$\textit{quasi-convex on} $X$ (see [9]),
iff for any $x_{1},x_{2}\in X$ and $\lambda \in \lbrack 0,1],$ either $%
F(x_{1})\subset F(\lambda x_{1}+(1-\lambda )x_{2})+S$ or $F(x_{2})\subset
F(\lambda x_{1}+(1-\lambda )x_{2})+S.$

(ii) $F$ is said to be (in the sense of [12 , Definition 3.6]) \textit{%
type-(v) properly }$S-$\textit{quasi-convex on} $X$ (see [9]), iff for any $%
x_{1},x_{2}\in X$ and $\lambda \in \lbrack 0,1],$ either $F(\lambda
x_{1}+(1-\lambda )x_{2})\subset F(x_{1})-S$ or $F(\lambda x_{1}+(1-\lambda
)x_{2})\subset F(x_{2})-S.$

If $-F$ is a type-(iii) [resp. type-(v)] $S-$properly quasiconvex set-valued
map, then, $F$ is said be type-(iii) [resp. type-(v)] $S-$properly
quasi-concave, which is equivalent to type-(iii) [resp. type-(v)] $(-S)$%
-properly quasi-convex set valued map.

(iii) $F:X\rightrightarrows Y$ is said to be (in the sense of [12 ,
Definition 3.6]) \textit{type-(iii)} \emph{naturally S-quasi-convex }on $X$,
iff for any $x_{1}$,$x_{2}\in X$ and $\lambda \in \lbrack 0,1],$ co$%
(F(x_{1})\cup F(x_{2}))\subset $ $F(\lambda x_{1}+(1-\lambda )x_{2})+S$.

iv) $F:X\rightrightarrows Y$ is said to be (in the sense of [12 , Definition
3.6]) \textit{type-(v)} \emph{naturally S-quasi-convex }on $X$, iff for any $%
x_{1}$,$x_{2}\in X$ and $\lambda \in \lbrack 0,1],$ $F(\lambda
x_{1}+(1-\lambda )x_{2})\subset $co$(F(x_{1})\cup F(x_{2}))-S$.

$F$ is said to be \textit{type-(iii)} \textit{[resp. type-(v)]} \textit{%
naturally }$S-$\textit{quasi-concave on} $X$, iff $-F$ is type-(iii) [resp.
type-(v)] naturally $S-$quasi-convex on $X.$

(v) $F:X\rightrightarrows Y$ is said to be \emph{S-quasi-convex }on $X$ (see
[24]), iff for any $x_{1}$,$x_{2}\in X$ and $\lambda \in \lbrack 0,1],$ ($%
F(x_{1})+S)\cap (F(x_{2})+S)\subset $ $F(\lambda x_{1}+(1-\lambda )x_{2})+S$.

(vi) $F$ is \textit{quasi-convex} $X$ [24] iff, for each $n$ and for every $%
x_{1},x_{2},...,x_{n}\in X,$ $\lambda =(\lambda _{1},\lambda
_{2},...,\lambda _{n})\in \Delta _{n-1},$ $\tbigcap\limits_{i=1}^{n}F(x_{i})%
\subset F(\tsum_{i=1}^{n}\lambda _{i}x_{i}).$

$F$ is said to be \textit{quasi-concave on} $X$, iff $-F$ is quasi-convex on 
$X.$

\textbf{Definition 2.3 }(see [26]) Let $X$ be a non-empty convex subset of a
topological vector space\textit{\ }$E$, let $Y$ be a subset of a topological
vector space $Z$ and $S$ a pointed closed convex cone in $Z$ with its
interior int$S\neq \emptyset .$ A vector-valued mapping $f:X\rightarrow Y$
is said to be \textit{natural}\emph{\ }$S-$\textit{quasi-convex} on $X$ iff $%
f(\lambda x_{1}+(1-\lambda )x_{2})\in $co$\{f(x_{1}),f(x_{2})\}-S$ for every 
$x_{1},x_{2}\in X$ and $\lambda \in \lbrack 0,1].$ This condition is
equivalent with the following condition: there exists $\mu \in \lbrack 0,1]$
such that $f(\lambda x_{1}+(1-\lambda )x_{2})\leq _{S}\mu f(x_{1})+(1-\mu
)f(x_{2}),$ where $x\leq _{S}y$ $\Leftrightarrow $ $y-x\in S.$

A vector-valued mapping f is said to be \emph{natural }$S-$\emph{%
quasi-concave} on $X$ if $-f$ is natural quasi $S-$convex on $X$.

\textit{Notation. }We will denote by $\Delta _{n-1}$ the standard
(n-1)-dimensional simplex in $R^{n},$ that is

$\Delta _{n-1}=\left\{ (\lambda _{1},\lambda _{2},...,\lambda _{n})\in R^{n}:%
\overset{n}{\underset{i=1}{\tsum }}\lambda _{i}=1\text{ and }\lambda
_{i}\geqslant 0,i=1,2,...,n\right\} .$

In this paper, we will also use the following notation:

$C^{\ast }(\Delta _{n-1})=\{g=(g_{1},g_{2},...,g_{n}):\Delta
_{n-1}\rightarrow \Delta _{n-1}$ where $g_{i}$ is continuous, $g_{i}(1)=1$
and $g_{i}(0)=0$ for each\textit{\ }$i\in \{1,2,...,n\}\}$

\textbf{Definition 2.4 }(see [3])\textit{\ }Let $X$ be a non-empty convex
subset of a topological vector space $E$ \ and $Y$ a non-empty subset of $E$.%
\textit{\ }The set-valued map\textit{\ }$F:X\rightrightarrows Y$ is said to
have \emph{weakly convex graph} (in short, it is a WCG correspondence) if,
for each $n\in N$ and for each finite set $\{x_{1},x_{2},...,x_{n}\}\subset
X $, there exists $y_{i}\in F(x_{i})$, $(i=1,2,...,n)$ such that

$\ \ \ \ \ \ \ \ \ \ \ \ \ \ \ \ \ \ \ \ \ \ \ \ \ \ \ (1.1)\ \ \ \ \ $co$%
(\{(x_{1},y_{1}),(x_{2},y_{2}),...,(x_{n},y_{n})\})\subset $Gr$(F)$

The relation (1.1) is equivalent to

$\ \ \ \ \ \ \ \ \ \ \ \ \ \ \ \ \ \ \ \ \ \ \ \ \ \ \ (1.2)\ \ \ \ \overset{%
n}{\underset{i=1}{\tsum }}\lambda _{i}y_{i}\in F(\overset{n}{\underset{i=1}{%
\tsum }}\lambda _{i}x_{i})\ \ \ \ \ \ \ (\forall (\lambda _{1},\lambda
_{2},...,\lambda _{n})\in \Delta _{n-1}).$

In [22] we introduced\ the concept of weakly naturally quasi-concave
set-valued map.

\textbf{Definition 2.5 }(see [22])Let $X$ be a non-empty convex subset of a
topological vector space $E$ \ and $Y$ a non-empty subset of a topological
vector space $Z$.\textit{\ }The set-valued map\textit{\ }$%
F:X\rightrightarrows Y$ is said to be \emph{weakly naturally quasi-concave
(WNQ) }iff, for each $n$ and for each finite set $\{x_{1},x_{2},...,x_{n}\}%
\subset X$, there exists $y_{i}\in F(x_{i})$, $(i=1,2,...,n)$ and $g\in
C^{\ast }(\Delta _{n-1})$ such that $\overset{n}{\underset{i=1}{\tsum }}%
g_{i}(\lambda _{i})y_{i}\in F(\overset{n}{\underset{i=1}{\tsum }}\lambda
_{i}x_{i})$ for every $(\lambda _{1},\lambda _{2},...,\lambda _{n})\in
\Delta _{n-1}.$

\textit{Remark 2.1} If $g_{i}(\lambda _{i})=\lambda _{i}$ for each $i\in
(1,2,...,n)$ and $(\lambda _{1},\lambda _{2},...,\lambda _{n})\in \Delta
_{n-1},$ we get a set-valued map with weakly convex graph, as it is defined
by Ding and He Yiran in [3]. In the same time, the weakly naturally
quasi-concavity is a weakening of the notion of naturally S-quasi-concavity
with $S=\{0\}.$

\textit{Remark 2.2} If $F$ is a single-valued mapping, then, it must be
natural $S$-quasiconcave for $S=\{0\}.$

\textit{Example 2.1} (see [22]) Let $F:[0,4]\rightrightarrows \lbrack -2,2]$
be defined by

$F(x)=\left\{ 
\begin{array}{c}
\lbrack 0,2]\text{ if }x\in \lbrack 0,2); \\ 
\lbrack -2,0]\text{ \ \ if \ }x=2; \\ 
(0,2]\text{ if }x\in (2,4].%
\end{array}%
\right. $

$F$ is neither upper semicontinuous, nor lower semicontinuous in $2.$ $F$
has not either got a weakly convex graph, since, if we consider $n=2,$ $%
x_{1}=1$ and $x_{2}=3,$ we have that co$\{(1,y_{1}),(3,y_{2})\}\nsubseteq $Gr%
$F,$ for every $y_{1}\in F(x_{1}),y_{2}\in F(x_{2}).$ We notice that $F$ is
not naturally $\{0\}-$quasi-concave, but it is weakly naturally
quasi-concave.

We proved in [22] the following fixed point theorem.

\textbf{Theorem 2.1} (see [22])\textit{Let }$Y$\textit{\ be a non-empty
subset of a topological vector space }$E$\textit{\ and }$K$\textit{\ be a }$%
(n-1)$\textit{- dimensional simplex in }$E$\textit{. Let }$%
F:K\rightrightarrows Y$\textit{\ be an weakly naturally quasi-concave
set-valued map and }$s:Y\rightarrow K$\textit{\ be a continuous function.
Then, there exists }$x^{\ast }\in K$\textit{\ such that }$x^{\ast }\in
s\circ F(x^{\ast })$\textit{.\medskip }

\textbf{3. Set-valued Maps with Generalized Convexity}\textit{\medskip }

In this section, we introduce several classes of cone convexity in order to
generalize the requirements for results concerning minimax inequalities.
Concerning the minimax problems we consider in this paper, we must underline
the behaviour importance of the set-valued maps $F(\cdot ,\cdot ):X\times
X\rightarrow Y$ in the points where their values contain or not maximal
(resp. minimal) elements of the certain sets of type $\tbigcup\limits_{y\in
X}F(x,y)$ or $\tbigcup\limits_{x\in X}F(x,y)$. We obtain the new definitions
through transferring the convexity properties of the maps from a variable to
another and by taking into consideration the maximal (resp. minimal)
elements. The reasons for our conception of generalized convex set-valued
maps come from the motivating work in the framework of minimax theory, where
the new properties prove to be necessary in order to obtain results by
giving up the continuity assumptions.

We firstly define the $S-$transfer $\mu -$convexity.\textit{\ }

\textbf{Definition 3.1 }Let $X$ be a convex set of a topological vector space%
\textit{\ }$E,$ let $Y$ be a non-empty set in the topological vector space $%
Z $ and let $F:X\times X\rightrightarrows Y$ be a set valued map with
non-empty values. $F$ is called $S-$\textit{transfer type-(v) }$\mu -$%
\textit{convex} \textit{in the first argument on} $X\times X$ iff, for each $%
n\in N$, $x_{1},x_{2},...,x_{n}\in X$ and $z\in X,$ we have that, for each $%
i\in \{1,2,...,n\},$ there exists $z_{i}=z_{i}(x_{1},x_{2},...,x_{n},z)\in X$
such that:\newline
i) $F(\tsum_{i=1}^{n}\lambda _{i}x_{i},z)\cap (\tbigcup_{y\in
X}F(x_{i},y))\subset F(x_{i},z_{i})-S$ for each $\lambda =(\lambda
_{1},\lambda _{2},...,\lambda _{n})\in \Delta _{n-1}$ with the property that 
$F(\tsum_{i=1}^{n}\lambda _{i}x_{i},z)\cap $Max$_{w}\tbigcup_{y\in
X}F(\tsum_{i=1}^{n}\lambda _{i}x_{i},y)\neq \emptyset $ or,

ii) $F(\tsum_{i=1}^{n}\lambda _{i}x_{i},z)\cap (\tbigcup_{y\in
X}F(x_{i},y))\subset F(x_{i},z_{i})-$int$S$ for each $\lambda =(\lambda
_{1},\lambda _{2},...,\lambda _{n})\in \Delta _{n-1}$ with the property that 
$F(\tsum_{i=1}^{n}\lambda _{i}x_{i},z)\cap $Max$_{w}\tbigcup_{y\in
X}F(\tsum_{i=1}^{n}\lambda _{i}x_{i},y)=\emptyset .$

$F$ is called $S-$\textit{transfer type-(v) }$\mu -$\textit{concave} \textit{%
in the first argument on} $X\times X$ if $-F$ is $S-$\textit{transfer
type-(v) }$\mu -$\textit{convex} \textit{in the first argument on} $X\times
X.\medskip $

\textit{Remark 3.1 }We can similarily define the $S-$transfer type-(iii)%
\textit{\ }$\mu -$convex set-valued maps.

\textit{Example 3.1} Let $X=[0,1],$ $Y=[-1,1],$ $S=[0,\infty )$ and $%
F:X\times X\rightrightarrows Y$ be defined by $F(x,y)=\left\{ 
\begin{array}{c}
\lbrack -1,y]\text{ if }0\leq x\leq y\leq 1; \\ 
\lbrack -x,y]\text{ if }0\leq y<x\leq 1.%
\end{array}%
\right. $

We will prove that $F$ is $S-$transfer type-(v)\textit{\ }$\mu -$convex in
the first argument$.$

Let $x_{1},x_{2},...,x_{n}\in X$ and $z\in Y.$ For each $i\in \{1,2,...,n\},$
$\tbigcup_{y\in X}F(x_{i},y)=[-1,1].$

Moreover, by computing, we obtain Max$_{w}\tbigcup_{y\in X}F(x_{i},y)=\{1\}$
and

$F(\tsum_{i=1}^{n}\lambda _{i}x_{i},z)=\left\{ 
\begin{array}{c}
\lbrack -1,z]\text{ \ \ \ \ \ \ \ if \ \ \ \ \ \ }0\leq
\tsum_{i=1}^{n}\lambda _{i}x_{i}\leq z\leq 1; \\ 
\lbrack -\tsum_{i=1}^{n}\lambda _{i}x_{i},z]\text{ if }0\leq
z<\tsum_{i=1}^{n}\lambda _{i}x_{i}\leq 1.%
\end{array}%
\right. $

For each $i\in \{1,2,...,n\},$ there exists $z_{i}\in Y,$ $z_{i}\geq \max
\{z,x_{i}\},$ so that $F(x_{i},z_{i})=[-1,z_{i}]$ and then:

i) if $z=1,$ $F(\tsum_{i=1}^{n}\lambda _{i}x_{i},z)\cap $Max$%
_{w}\tbigcup_{y\in X}F(x_{i},y)=F(\tsum_{i=1}^{n}\lambda _{i}x_{i},z)\cap
\{1\}\neq \emptyset $ and $F(\tsum_{i=1}^{n}\lambda _{i}x_{i},z)\subset
F(x_{i},z_{i})-S$ or

ii) if $z<1,$ $F(\tsum_{i=1}^{n}\lambda _{i}x_{i},z)\cap $Max$%
_{w}\tbigcup_{y\in X}F(x_{i},y)=F(\tsum_{i=1}^{n}\lambda _{i}x_{i},z)\cap
\{1\}=\emptyset $ and $F(\tsum_{i=1}^{n}\lambda _{i}x_{i},z)\subset
F(x_{i},z_{i})-$int$S$ $\medskip $

\noindent \textit{Remark 3.2} The $S-$transfer type-(v) $\mu -$convexity in
the first argument is implied by the following property, which we call $%
\alpha :$

$(\alpha ):$ For each $x\in X,$ $A_{x}=\cup _{y\in X}F(x,y)$ is compact and
there exists $z_{x}\in Z$ such that $z_{x}\in $Max$\cup _{y\in X}F(x,y)$ and 
$\cup _{y\in X}F(x,y)\subset z_{x}-S.$

We note that according to Lemma 2.1, $\cup _{y\in X}F(x,y)\subset $Max$\cup
_{y\in X}F(x,y)-S.$

The $S-$transfer type-(v) $\mu -$concavity in the second argument is implied
by the following property $\alpha ^{\prime }:$

$(\alpha ^{\prime }):$ For each $y\in X,$ $A_{y}=\cup _{x\in X}F(x,y)$ is
compact and there exists $z_{y}\in Z$ such that $z_{y}\in $Max$\cup _{x\in
X}F(x,y)$ and $\cup _{x\in X}F(x,y)\subset z_{y}+S.$

The set valued map from Example 3.1 verifies the property $\alpha $.

The condition $\alpha $ is not fulfilled in the next example.

\textit{Example 3.2} Let $S((0,0),x)=\{(u,v)\in \lbrack -1,1]\times \lbrack
-1,1]:u^{2}+v^{2}\leq x^{2}\},$

$S_{+}((0,0),x)=\{(u,v)\in \lbrack 0,1]\times \lbrack -1,1]:u^{2}+v^{2}\leq
x^{2}\}$ and

$S_{-}((0,0),x)=\{(u,v)\in \lbrack -1,0]\times \lbrack -1,1]:u^{2}+v^{2}\leq
x^{2}\}.$

Let us define $F:[0,1]\times \lbrack 0,1]\rightrightarrows \lbrack
-1,1]\times \lbrack -1,1]$ by

$F(x,y)=\left\{ 
\begin{array}{c}
S((0,0),1)\text{ \ \ \ \ if \ \ \ \ }x=1\text{ \ \ \ \ and }y\in \lbrack
0,1]. \\ 
S_{+}((0,0),x)\text{ \ if }0<x<1\text{ \ and \ }x\leq y\leq 1; \\ 
S_{-}((0,0),x)\text{ \ \ \ \ \ \ \ \ \ \ if \ \ \ \ \ \ \ \ \ \ \ }0<y<x<1;
\\ 
\{(0,0)\}\text{ \ \ \ \ \ \ if }x=0\text{ \ \ \ \ \ and \ \ \ \ \ }y\in
\lbrack 0,1].%
\end{array}%
\right. $

$F$ is $R_{+}^{2}$-transfer type-(v)\textit{\ }$\mu $ convex in the first
argument.$\medskip $

The Definition 3.1 can be weakened in the following way.

\textbf{Definition 3.2 }Let $X$ be a convex set of a topological vector space%
\textit{\ }$E,$ let $Y$ be a non-empty set in the topological vector space $%
Z $ and let $F:X\times X\rightrightarrows Y$ be a set-valued map with
non-empty values. $F$ is called $S-$\textit{transfer weakly type-(v) }$\mu -$%
\textit{convex} \textit{in the first argument on} $X\times X$ iff, for each $%
n\in N$, $x_{1},x_{2},...,x_{n}\in X$ and $z\in X,$ we have that, there
exist $i_{0}\in \{1,2,...,n\}$ and $%
z_{i_{0}}=z_{i_{0}}(x_{1},x_{2},...,x_{n},z)\in X$ such that:\newline
i) $F(\tsum_{i=1}^{n}\lambda _{i}x_{i},z)\cap (\tbigcup_{y\in
X}F(x_{i_{0}},y))\subset F(x_{i_{0}},z_{i_{0}})-S$ for each $\lambda
=(\lambda _{1},\lambda _{2},...,\lambda _{n})\in \Delta _{n-1}$ with the
property that $F(\tsum_{i=1}^{n}\lambda _{i}x_{i},z)\cap $Max$%
_{w}\tbigcup_{y\in X}F(\tsum_{i=1}^{n}\lambda _{i}x_{i},y)\neq \emptyset $
or,

ii) $F(\tsum_{i=1}^{n}\lambda _{i}x_{i},z)\cap (\tbigcup_{y\in
X}F(x_{i_{0}},y))\subset F(x_{i_{0}},z_{i_{0}})-$int$S$ for each $\lambda
=(\lambda _{1},\lambda _{2},...,\lambda _{n})\in \Delta _{n-1}$ with the
property that $F(\tsum_{i=1}^{n}\lambda _{i}x_{i},z)\cap $Max$%
_{w}\tbigcup_{y\in X}F(\tsum_{i=1}^{n}\lambda _{i}x_{i},y)=\emptyset .$

$F$ is called $S-$\textit{transfer weakly type-(v) }$\mu -$\textit{concave} 
\textit{in the first argument on} $X\times X$ if $-F$ is $S-$\textit{%
transfer weakly type-(v) }$\mu -$\textit{convex} \textit{in the first
argument on} $X\times X.\medskip $

\textit{Remark 3.3.}We can\textit{\ }similarily\textit{\ }define the\textit{%
\ }$S$-transfer weakly type-(iii)\textit{\ }$\mu $ convex set-valued maps.

\textit{Remark 3.4. }If $F:X\times X\rightarrow Z$ is type-(v) properly $S-$%
quasi-convex in the first argument, then, $F$ is $S$-transfer weakly type-(v)%
\textit{\ }$\mu $ convex in the first argument.

Indeed, let $x_{1},x_{2},...,x_{n}\in X$ and $\lambda =(\lambda _{1},\lambda
_{2},...,\lambda _{n})\in \Delta _{n-1}$. We have that $F(\tsum%
\limits_{i=1}^{n}\lambda _{i}x_{i},y)\subset F(x_{i_{0}},y)-S$ for each $%
\lambda \in \Delta _{n-1}$, $y\in X$ and an idex $i_{0}\in \{1,2,...,n\}.$
Then, for each $z\in X,$ there exists $z_{i_{0}}=z$ such that $%
F(\tsum\limits_{i=1}^{n}\lambda _{i}x_{i},z)\cap (\tbigcup\limits_{z\in
X}F(x_{i_{0}},z))\subset F(x_{i_{0}},z_{i_{0}})-S.$

Consequently, the notion of $S-$transfer weakly type-(v)\textit{\ }$\mu -$%
convexity is weaker than the type-(v) properly $S-$quasi-convexity and, in
certain cases, it is implied by the property $\alpha .$

\textit{Example 3.3} Let $X=[0,1],$ $Y=[-1,1],$ $S=[0,\infty )$ and $%
F:X\times X\rightrightarrows Y$ be defined by $F(x,y)=\left\{ 
\begin{array}{c}
\lbrack 0,y]\text{ if }0\leq x\leq y\leq 1; \\ 
\lbrack -x,y]\text{ if }0\leq y<x\leq 1.%
\end{array}%
\right. $

$F$ is $S-$transfer weakly type-(v)\textit{\ }$\mu -$convex in the first
argument$.$

Now, we are introducing a similar definition for single valued mappings.

\textbf{Definition 3.3 }Let $X$ be a convex set of a topological vector space%
\textit{\ }$E$ and let $Y$ be a non-empty set in the topological vector
space $Z.$

The mapping $f:X\times X\rightarrow Y$ is called $S-$\textit{transfer }$\mu
- $\textit{convex} \textit{in the first argument on} $X\times X$ iff, for
each $n\in N$, $x_{1},x_{2},...,x_{n}\in X$ and $z\in X,$ we have that, for
each $i\in \{1,2,...,n\},$ there exists $%
z_{i}=z_{i}(x_{1},x_{2},...,x_{n},z)\in X$ such that, if $%
f(\tsum_{i=1}^{n}\lambda _{i}x_{i},z)\in \tbigcup_{y\in X}f(x_{i},y)$ for
each $\lambda =(\lambda _{1},\lambda _{2},...,\lambda _{n})\in \Delta _{n-1},
$ the following condition is fulfilled:

i) $f(\tsum_{i=1}^{n}\lambda _{i}x_{i},z)\in f(x_{i},z_{i})-S$ for each $%
\lambda =(\lambda _{1},\lambda _{2},...,\lambda _{n})\in \Delta _{n-1}$ with
the property that $f(\tsum_{i=1}^{n}\lambda _{i}x_{i},z)\in $Max$%
_{w}(\tbigcup_{y\in X}f(\tsum_{i=1}^{n}\lambda _{i}x_{i},y))$ or,

ii) $f(\tsum_{i=1}^{n}\lambda _{i}x_{i},z)\in f(x_{i},z_{i})-$int$S$ for
each $\lambda =(\lambda _{1},\lambda _{2},...,\lambda _{n})\in \Delta _{n-1}$
with the property that $f(\tsum_{i=1}^{n}\lambda _{i}x_{i},z)\notin $Max$%
_{w}(\tbigcup_{y\in X}F(\tsum_{i=1}^{n}\lambda _{i}x_{i},y)).$

The mapping $f$ is called $S-$\textit{transfer }$\mu -$\textit{concave} 
\textit{in the first argument on} $X\times X$ iff $-f$ is $S-$\textit{%
transfer }$\mu -$\textit{convex} \textit{in the first argument on} $X\times
X.$

\textit{Example 3.4} Let $X=[0,1],$ $Y=[-1,0],$ $S=[0,\infty )$ and $%
f:X\times X\rightarrow Y$ be defined by $f(x,y)=\left\{ 
\begin{array}{c}
1\text{ if }0\leq x\leq y\leq 1; \\ 
x\text{ if }0\leq y<x\leq 1.%
\end{array}%
\right. $

We will prove that $f$ is $S-$transfer $\mu -$convex in the first argument$.$

Let $x_{1},x_{2},...,x_{n},z\in X.$ For each $i\in \{1,2,...,n\},$ $%
\tbigcup_{y\in X}f(x_{i},y)=\{x_{i},1\},$ Max$_{w}$ $\tbigcup_{y\in
X}f(x_{i},y)=\{1\}$ and we have that, for each $\lambda =(\lambda
_{1},\lambda _{2},...,\lambda _{n})\in \Delta _{n-1},$ if $%
f(\tsum_{i=1}^{n}\lambda _{i}x_{i},z)\in \{x_{i},1\},$ there exists $%
z_{i}\in Y,$ $z_{i}\geq \max \{z,x_{i}\},$ so that $f(x_{i},z_{i})=1,$ and
then:

i) if $z=1$ and $f(\tsum_{i=1}^{n}\lambda _{i}x_{i},z)=1$ for $(\lambda
_{1},\lambda _{2},...,\lambda _{n})\in \Delta _{n-1},$ we have that $%
f(\tsum_{i=1}^{n}\lambda _{i}x_{i},z)\in f(x_{i},z_{i})-S$ or,

ii) if $z<1$ and $f(\tsum_{i=1}^{n}\lambda _{i}x_{i},z)\neq 1$ for $(\lambda
_{1},\lambda _{2},...,\lambda _{n})\in \Delta _{n-1},$ we have that $%
f(\tsum_{i=1}^{n}\lambda _{i}x_{i},z)\in f(x_{i},z_{i})-$int$S.\medskip $

The next notion is stronger than the properly $S-$quasi-convexity and it is
adapted for set-valued maps with two variables. We consider pairs of points
in the product space $X\times X$. We keep constant one component and we
consider any convex combination of the other ones. By comparing the images
of $F$ in all these pairs of points, we obtain the following definition.

\textbf{Definition 3.4} Let $X$ be a non-empty convex subset of a
topological vector space\textit{\ }$E,$ $Z$ a real topological vector space
and $S$ a pointed closed convex cone in $Z$ with its interior int$S\neq
\emptyset .$ Let $F:X\rightrightarrows Z$ be a set-valued map with non-empty
values.

(i) $F$ is said to be \textit{type-(iii) pair} \textit{properly }$S-$\textit{%
quasi-convex on} $X\times X$ \textit{in the first argument,} iff, for any $%
(x_{1},y_{1}),(x_{2},y_{2})\in X\times X$ and $\lambda \in \lbrack 0,1],$
either $F(x_{1},y_{1})\subset F(\lambda x_{1}+(1-\lambda )x_{2},y_{1})+S$ or 
$F(x_{2},y_{2})\subset F(\lambda x_{1}+(1-\lambda )x_{2},y_{2})+S.$

(ii) $F$ is said to be \textit{type-(v) pair} \textit{properly }$S-$\textit{%
quasi-convex on} $X\times X$ \textit{in the first argument,} iff, for any $%
(x_{1},y_{1}),(x_{2},y_{2})\in X\times X$ and $\lambda \in \lbrack 0,1],$
either $F(\lambda x_{1}+(1-\lambda )x_{2},y_{1})\subset F(x_{1},y_{1})-S$ or 
$F(\lambda x_{1}+(1-\lambda )x_{2},y_{2})\subset F(x_{2},y_{2})-S.$

(iii) $F$ is said to be \textit{type-(iii) [resp. type-(v)] pair} \textit{%
properly }$S-$\textit{quasi-concave on} $X$ in the first argument, iff, $-F$
is type-(iii) [resp. type-(v)] pair\textit{\ }properly $S-$quasi-convex in
the first argument\textit{\ }on $X.$

(iv) $F$ is said to be \textit{pair} \textit{properly quasi-convex }iff for
any $(x_{1},y_{1}),(x_{2},y_{2})\in X\times X$ and $\lambda \in \lbrack
0,1], $ either $F(x_{1},y_{1})\subset F(\lambda x_{1}+(1-\lambda
)x_{2},y_{1})$ or $F(x_{2},y_{2})\subset F(\lambda x_{1}+(1-\lambda
)x_{2},y_{2}).$

$\mathit{F}$ is said to be \textit{pair} \textit{properly quasi-concave if }$%
-F$ is pair properly $S-$quasi-convex.$\medskip $

\textit{Example 3.5} Let $X=[0,1],$ $Y=[-1,1],$ $S=[0,\infty )$ and $%
F:X\times X\rightrightarrows Y$ be defined by $F(x,y)=\left\{ 
\begin{array}{c}
\lbrack -1,1]\text{ if }0\leq x\leq y\leq 1; \\ 
\lbrack -x,1]\text{ if }0\leq y<x\leq 1.%
\end{array}%
\right. $

$F$\ is type-(iii)\textit{\ }pair properly quasi-concave in the second
argument on $X.$

\textit{Remark 3.5. }$S-$transfer $\mu -$convexity does not imply pair
properly $S-$quasi-convexity. The set valued map from Example 3.2 is $%
R_{+}^{2}-$transfer type-(v)\textit{\ }$\mu -$convex in the first argument,
but it is not type-(v)\textit{\ }pair properly $R_{+}^{2}-$quasi-convex in
the first argument.

If we consider $(x_{1},y_{1})=(\frac{1}{15},\frac{9}{10}),$ $(x_{2},y_{2})=(%
\frac{1}{4},\frac{1}{5})$ and $x_{0}=\frac{1}{5}\in $co$\{x_{1},x_{2}\},$
then, $F(x_{1},y_{1})=S_{+}((0,0),\frac{1}{15}),$ $%
F(x_{2},y_{2})=S_{-}((0,0),\frac{1}{4}),$ $F(x_{0},y_{1})=S_{+}((0,0),\frac{1%
}{5})$ and $F(x_{0},y_{2})=S_{+}((0,0),\frac{1}{5}).$ It follows that
neither $F(x_{0},y_{1})\subset F(x_{1},y_{1})-R_{+}^{2},$ nor $%
F(x_{0},y_{2})\subset F(x_{2},y_{2})-R_{+}^{2}$ and then, $F$ is not
type-(v) pair properly $R_{+}^{2}-$quasi-convex in the first argument.

Conversely, the pair properly $S-$quasi-convexity does not imply $S-$%
transfer $\mu -$convexity. The following example is concludent in this
respect.

\textit{Example 3.6.} For each $(x,y)\in \lbrack 0,1]\times \lbrack 0,1],$
let us define

$S((0,y),x)=\{(u,v)\in R^{2}\times R^{2}:u^{2}+(v-y)^{2}\leq x^{2}\}$ and

$S((y,0),x)=\{(u,v)\in \in R^{2}\times R^{2}:(u-y)^{2}+v^{2}\leq x^{2}\}.$

Let $S=R_{+}^{2}$ and $F:[0,1]\times \lbrack 0,1]\rightrightarrows \lbrack
-2,2]\times \lbrack -2,2]$ be defined by

$F(x,y)=\left\{ 
\begin{array}{c}
S((0,y),x)\text{ \ \ \ \ if \ \ }(x,y)\in \lbrack 0,1]\times ([0,1]\cap Q);
\\ 
S((y,0),x)\text{ if }(x,y)\in \lbrack 0,1]\times ([0,1]\cap (R\backslash Q)).%
\end{array}%
\right. $

The set valued map $F$ is type-(v) pair properly $R_{+}^{2}-$quasi-convex in
the first argument, but it is not $R_{+}^{2}-$transfer type-(v) $\mu -$%
convex in the first argument.

Indeed, let us consider first $(x_{1},y_{1})$ and $(x_{2},y_{2})\in \lbrack
0,1].$ Without loss of generalization, we can assume that $x_{1}\leq
x(\lambda )\leq x_{2}$ for each $\lambda \in \lbrack 0,1]$, where $x(\lambda
)=\lambda x_{1}+(1-\lambda )x_{2}.$ Consequently, $F(x(\lambda
),y_{2})\subset F(x_{2},y_{2})-S$ and $F$ is type-(v) pair properly $%
R_{+}^{2}-$quasi-convex in the first argument.

In order to prove the second assertion, let us consider $x_{1},x_{2}\in
\lbrack 0,1]$ and

$x(\lambda )=\lambda x_{1}+(1-\lambda )x_{2},$ where $\lambda \in \lbrack
0,1].$

For $i=1,2$ and $y=0$ the following equality holds: $F(x(\lambda ),0)\cap
\tbigcup\limits_{y\in \lbrack 0,1]}F(x_{i},y)=$

$(F(x(\lambda ),0)\cap \tbigcup\limits_{y\in \lbrack 0,1]\cap
Q}F(x_{i},y))\cup (F(x(\lambda ),0)\cap \tbigcup\limits_{y\in \lbrack
0,1]\cap (R\backslash Q)}F(x_{i},y))$

and there is not any $z_{i}\in \lbrack 0,1]$ such that $F(x(\lambda ),0)\cap
\tbigcup\limits_{y\in \lbrack 0,1]}F(x_{i},y)\subset
F(x_{i},z_{i})-R_{+}^{2}.$

We conclude that $F$ is not $R_{+}^{2}-$transfer type-(v) $\mu -$convex in
the first argument.

For single valued mappings, the next definition is proposed.

\textbf{Definition 3.5} Let $X$ be a nonempty convex subset of a topological
vector space\textit{\ }$E,$ $Z$ a real topological vector space and $S$ a
pointed closed convex cone in $Z$ with its interior int$S\neq \emptyset .$
Let $f:X\rightarrow Z$ be a set-valued map with non-empty values.

(i) $f$ is said to be \textit{pair} \textit{properly }$S-$\textit{%
quasi-convex on} $X\times X$ \textit{in the first argument,} iff, for any $%
(x_{1},y_{1}),(x_{2},y_{2})\in X\times X$ and $\lambda \in \lbrack 0,1],$
either $f(x_{1},y_{1})\subset f(\lambda x_{1}+(1-\lambda )x_{2},y_{1})+S$ or 
$f(x_{2},y_{2})\subset f(\lambda x_{1}+(1-\lambda )x_{2},y_{2})+S.$

$f$ is said to be \textit{pair} \textit{properly }$S-$\textit{quasi-concave
in the first argument on} $X\times X$\textit{,} iff $-f$ is properly $S-$%
quasi-convex in the first argument\textit{\ }on $X\times X.$

The usual naturally $S-$quasi-convexity requirement in the minimax
inequalities for set-valued maps can be weakened. In the definition we
propose below, we take into consideration the bahaviour of the set-valued
maps in the points where their values do not contain minimal (resp. maximal)
points of some certain sets of $\tbigcup_{y\in X}F(x,y)$ or $\tbigcup_{x\in
X}F(x,y)$ types$.$

\textbf{Definition 3.6 }Let $X$ be a convex set of a topological vector space%
\textit{\ }$E,$ let $Y$ be a non-empty set in the topological vector space $%
Z $ and let $F:X\times X\rightrightarrows Y$ be a set-valued map with
non-empty values.

i) $F$ is called \textit{transfer type-(iii) properly }$S-$\textit{%
quasi-convex} \textit{in the first argument} on $X\times X$ iff, for each
elements $x_{1},x_{2},z\in X,$ $\lambda \in (0,1)$ and $i\in \{1,2\}$, the
following condition is fulfilled: $F(\lambda x_{1}+(1-\lambda )x_{2},z)\cap $%
Min$_{w}(\tbigcup_{y\in X}F(x_{i},y))=\emptyset $ implies that $%
F(x_{i},z)\subset F(\lambda x_{1}+(1-\lambda )x_{2},z)+S$.

ii) $F$ is called \textit{transfer type-(v) properly }$S-$\textit{%
quasi-convex} \textit{in the first argument} on $X\times X$ iff, for each
elements $x_{1},x_{2},z\in X,$ $\lambda \in (0,1)$ and $i\in \{1,2\}$, the
following condition is fulfilled: $F(\lambda x_{1}+(1-\lambda )x_{2},z)\cap $%
Min$_{w}(\tbigcup_{y\in X}F(x_{i},y))=\emptyset $ implies $F(\lambda
x_{1}+(1-\lambda )x_{2},z)\subset F(x_{i},z)-S$.

$F$ is called \textit{transfer type-(iii) [resp.type-(v)] properly }$S-$%
\textit{quasi}-\textit{concave} \textit{in the first argument }on $X\times X$
if $-F$ is transfer type-(iii) [resp.type-(v)]\textit{\ }properly $S-$%
quasi-convex in the first argument\textit{\ }on $X\times X.$

\textit{Remark 3.6.} If $F(\cdot ,y)$ is naturally $S-$quasi-convex for each 
$y\in X,$ then, $F$ is transfer properly $S-$quasi-convex in the first
argument on $X\times X.$

\textit{Remark 7. }If $F$ is transfer properly $S-$quasi-convex in the first
argument on $X\times X,$ then, $F$ is $S-$transfer weakly $\mu -$convex in
the first argument.

Conversely, it is not true. The set-valued map $F$ defined in Example 3.2 is 
$S-$transfer weakly (type-v) $\mu -$convex in the first argument, but it is
not type-(v) transfer properly $S-$quasi-convex.

\textit{Example 3.7} Let $X=[0,1],$ $Y=[-1,1],$ $S=[0,\infty )$ and $%
F:X\times X\rightrightarrows Y$ be defined by $F(x,y)=\left\{ 
\begin{array}{c}
\lbrack 0,y]\text{ if }0\leq x\leq y\leq 1; \\ 
\lbrack -x,y]\text{ if }0\leq y<x\leq 1.%
\end{array}%
\right. $

We prove that\ $F(\cdot ,y)$\ is type-(iii)\textit{\ }naturally $S-$%
quasi-concave on $X$ (and then, $F$ is transfer type-(iii)\textit{\ }%
properly $S-$quasi-concave in the first argument on $X\times X$)$.$

Let $y\in \lbrack 0,1]$ be fixed, $x_{1},x_{2}\in \lbrack 0,1]$, $\lambda
\in \lbrack 0,1]$ and $x(\lambda )=\lambda x_{1}+(1-\lambda )x_{2}$.

1) If $x_{1}\geq x_{2}\geq y,$ then, $F(x_{1},y)=[-x_{1},y],$ $%
F(x_{2},y)=[-x_{2},y]$, $F(x(\lambda ),y)=[-x(\lambda ),y]$ and

co\{$F(x_{1},y),F(x_{2},y)\}=[-x_{1},y]\subset \lbrack -x(\lambda
),y]-[0,\infty )=F(x(\lambda ),y)-[0,\infty );$

2) if $x_{1}\leq x_{2}\leq y,$ then, $F(x_{1},y)=[0,y],$ $F(x_{2},y)=[0,y]$, 
$F(x(\lambda ),y)=[0,y]$ and

co\{$F(x_{1},y),F(x_{2},y)\}=[0,y]\subset \lbrack 0,y]-[0,\infty
)=F(x(\lambda ),y)-[0,\infty );$

3) if $x_{1}\geq y\geq x_{2},$ then, $F(x_{1},y)=[-x_{1},y]$, $%
F(x_{2},y)=[0,y]$ and

co\{$F(x_{1},y),F(x_{2},y)\}=[-x_{1},y];$

if $x_{1}\geq x(\lambda )\geq y\geq x_{2},$ then, $F(x(\lambda
),y)=[-x(\lambda ),y]$ and

co\{$F(x_{1},y),F(x_{2},y)\}=[-x_{1},y]\subset \lbrack -x(\lambda
),y]-[0,\infty )=F(x(\lambda ),y)-[0,\infty );$

if $x_{1}\geq y\geq x(\lambda )\geq x_{2},$ then, $F(x(\lambda ),y)=[0,y]$
and

co\{$F(x_{1},y),F(x_{2},y)\}=[-x_{1},y]\subset \lbrack 0,y]-[0,\infty
)=F(x(\lambda ),y)-[0,\infty ).$\medskip

The usual properly $S-$quasi-convexity assumption in the minimax theorems
with set-valued maps can be also generalized. In order to obtain necessary
conditions in our results, we introduce the following definitions.

\textbf{Definition 3.7 }Let $X$ be a convex set of a topological vector space%
\textit{\ }$E,$ let $Y$ be a non-empty set in the topological vector space $%
Z $ and let $F:X\times X\rightrightarrows Y$ be a set-valued map with
non-empty values. $F$ \textit{satisfies the condition }$\gamma $\textit{\ on}
$X\times X$ iff: \newline
$(\gamma )$ \ \ there exist $n\in N,$ $%
(x_{1},y_{1}),(x_{2},y_{2}),...,(x_{n},y_{n})\in X\times X$, $y^{\ast }\in $%
co$\{x_{1},x_{2},...,x_{n}\}$ such that $F(x_{i},y_{i})\subset
F(x_{i},y^{\ast })-S$ and $F(x_{i},y_{i})\cap $Max$_{w}\cup _{z\in
X}F(x_{i},z)\neq \emptyset $ for each $i\in \{1,2,...,n\}$.

\textit{Example 3.8} Let $X=[0,1],$ $Y=[-1,1],$ $S=[0,\infty )$ and $%
F:X\times X\rightrightarrows Y$ be defined by $F(x,y)=\left\{ 
\begin{array}{c}
\lbrack 0,y]\text{ if }0\leq x\leq y\leq 1; \\ 
\lbrack -x,y]\text{ if }0\leq y<x\leq 1.%
\end{array}%
\right. $

We prove that\ $F(x,\cdot )$\ satisfies the condition $\gamma .$ In fact,
there exist $(x_{1},y_{1})=(0,1),$ $(x_{2},y_{2})=(1,1)\in X\times X$ such
that $F(x_{i},y_{i})\cap $Max$_{w}\cup _{z\in X}F(x_{i},z)\neq \emptyset ,$ $%
i=1,2.$ There also exists $y^{\ast }=1\in $co$\{x_{1},x_{2}\}$ such that $%
[0,1]=F(x_{1},y_{1})\subset F(x_{1},y^{\ast })-[0,\infty )$ and $%
[0,1]=F(x_{2},y_{2})\subset F(x_{2},y^{\ast })-[0,\infty ).$

\textbf{Definition 3.8 }Let $X$ be a convex set of a topological vector space%
\textit{\ }$E,$ let $Y$ be a non-empty set in the topological vector space $%
Z $ and let $F:X\times X\rightrightarrows Y$ be a set-valued map with
non-empty values. $F$ \textit{satisfies the condition }$\gamma ^{\prime }$%
\textit{\ on} $X\times X$ iff:

$(\gamma ^{\prime })$ \ \ there exist $n\in N,$ $%
(x_{1},y_{1}),(x_{2},y_{2}),...,(x_{n},y_{n})\in X\times X$ and $x^{\ast
}\in $co$\{y_{1},y_{2},...,y_{n}\}$ such that $F(x_{i},y_{i})\subset
F(x_{i},y^{\ast })+S$ and $F(x_{i},y_{i})\cap $Min$_{w}\cup _{x\in
X}F(x,y_{i})\neq \emptyset $ for each $i\in \{1,2,...,n\}$.

\textbf{4. Minimax Theorems for Set-valued Maps without Continuity\medskip }

In this section, we establish some generalized Ky Fan minimax inequalities.

Firstly, we are proving the following lemma, which is comparable with Lemma
3.1 in [32], but our result does not involve continuity assumptions.
Instead, we use several generalized convexity properties for set-valued maps
introduced in Section 3.

Lemma 4.1 will be used to prove the minimax Theorem 4.1.\medskip

\textbf{Lemma 4.1} \textit{Let }$X$\textit{\ be a (n-1) dimensional simplex
of a Hausdorff topological vector space }$E,$\textit{\ }$Y$\textit{\ a
compact set in the Hausdorff topological vector space }$Z$\textit{\ and let }%
$S$\textit{\ be a pointed closed convex cone in }$Z$\textit{\ with its
interior int}$S\neq \emptyset .$\textit{\ Let }$F:X\times X\rightrightarrows
Y$\textit{\ be a set-valued map with non-empty values.}

\textit{(i) Let us suppose that }$\tbigcup_{y\in X}F(x,y)$\textit{\ is a
compact set for each }$x\in X$. \textit{If }$F$ \textit{is }$S-$\textit{%
transfer type-(v) }$\mu -$\textit{convex in the first argument on} $X\times
X $\textit{,} $F$\textit{\ is type-(iii) pair properly quasi-concave in the
second argument on }$X\times X$ \textit{and }$F(\cdot ,y)$\textit{\ is
type-(iii) naturally }$S-$\textit{quasi-concave on }$X$ \textit{for each }$%
y\in X,$\textit{\ then, there exists }$x^{\ast }\in X$\textit{\ such that }$%
F(x^{\ast },x^{\ast })\cap $\noindent Max$_{w}\tbigcup_{y\in X}F(x^{\ast
},y)\neq \emptyset .$

\textit{(ii) Suppose that }$\tbigcup_{x\in X}F(x,y)$\textit{\ is a compact
set for each }$y\in X$. \textit{If }$F$ \textit{is transfer type-(v) }$\mu -$%
\textit{concave in the second argument on} $X\times X,$ $F$\textit{\ is
type-(iii) pair properly quasi-convex in the first argument on }$X\mathit{%
\times X}$ \textit{and }$F(x,\cdot )$\textit{\ is type-(iii) naturally }$S-$%
\textit{quasi-convex\ on }$X$ \textit{for each }$x\in X,$\textit{\ then,
there exists }$y^{\ast }\in X$\textit{\ such that }$F(y^{\ast },y^{\ast
})\cap $Min$_{w}\cup _{x\in X}F(x,y^{\ast })\neq \emptyset .\medskip $

\textit{Proof.} (i) Let us define the set-valued map $T:X\rightrightarrows X$
by

$T(x)=\{y\in X:F(x,y)\cap $Max$_{w}\cup _{z\in X}F(x,z)\neq \emptyset \}$
for each $x\in X.$

We claim that $T$ is non-empty valued. Indeed, since $\cup _{z\in X}F(x,z)$
is a compact set for each $x\in X,$ according to Lemma 2.1, Max$_{w}\cup
_{z\in X}F(x,z)\neq \emptyset .$ For each $x\in X,$ let $z_{x}\in $Max$%
_{w}\cup _{z\in X}F(x,z).$ Then, there exists $y_{x}\in X$ such that $%
z_{x}\in F(x,y_{x}).$ It is clear that $y_{x}\in T(x)=\{y\in X:F(x,y)\cap $%
Max$_{w}\cup _{z\in X}F(x,z)\}$ and, consequently, $T(x)\neq \emptyset $ for
each $x\in X.$

Further,we will prove that $T$ is weakly naturally quasi-concave.

Let $x_{1},x_{2},...,x_{n}\in X.$ For each $i\in 1,...,n,$ there exists $%
y_{i}\in T(x_{i})$, that is $F(x_{i},y_{i})\cap $Max$_{w}\tbigcup_{z\in
X}F(x_{i},z)\neq \emptyset .$

By contrary, we assume that $T$ is not weakly naturally quasi-concave. Then,
for each $g\in C^{\ast }(\Delta _{n-1})$, there exists $\lambda
^{g}=(\lambda _{1}^{g},\lambda _{2}^{g},...,\lambda _{n}^{g})\in \Delta
_{n-1}$ such that $\tsum_{i=1}^{n}g_{i}(\lambda _{i}^{g})y_{i}\notin
T(\tsum_{i=1}^{n}\lambda _{i}^{g}x_{i}),$ relation which is equivalent with
the following one:

$F(\tsum_{i=1}^{n}\lambda _{i}^{g}x_{i},\tsum_{i=1}^{n}g_{i}(\lambda
_{i}^{g})y_{i})\cap $Max$_{w}\cup _{z\in X}F(\tsum_{i=1}^{n}\lambda
_{i}^{g}x_{i},z)=\emptyset .$

Since the set-valued map $F$ is $S-$transfer type-(v)\textit{\ }$\mu -$%
convex in the first argument and $F(\tsum_{i=1}^{n}\lambda
_{i}^{g}x_{i},\tsum_{i=1}^{n}g_{i}(\lambda _{i}^{g})y_{i})\cap $Max$%
_{w}F(\tsum_{i=1}^{n}\lambda _{i}^{g}x_{i},X)=\emptyset $, it follows that,
for each $i\in \{1,2,...,n\},$ there exists the element $z_{i_{0}}\in X$
such that the following relation is fulfilled: $F(\tsum_{i=1}^{n}\lambda
_{i}^{g}x_{i},\tsum_{i=1}^{n}g_{i}(\lambda _{i}^{g})y_{i})\cap
(\tbigcup_{z\in X}F(x_{i},z))\subset F(x_{i},z_{i_{0}})-$int$S.$

Let $t_{i}\in F(\tsum_{i=1}^{n}\lambda
_{i}^{g}x_{i},\tsum_{i=1}^{n}g_{i}(\lambda _{i}^{g})y_{i})\cap
(\tbigcup_{z\in X}F(x_{i},z))$ and $u_{i}\in F(x_{i},z_{i_{o}})$ such that $%
t_{i}=u_{i}-s_{i},$ $s_{i}\in $int$S.$ It follows that $u_{i}\in
\tbigcup_{z\in X}F(x_{i},z)\cap \{t_{i}+$int$S\}\neq \emptyset ,$ that is $%
t_{i}\in F(\tsum_{i=1}^{n}\lambda _{i}^{g}x_{i},\tsum_{i=1}^{n}g_{i}(\lambda
_{i}^{g})y_{i})\cap (\tbigcup_{z\in X}F(x_{i},z))$ implies the fact that $%
t_{i}\notin $Max$_{w}\cup _{z\in X}F(x_{i},z).$ Consequently, we have that,
for each index $i\in \{1,2,...,n\},$

$F(\tsum_{i=1}^{n}\lambda _{i}^{g}x_{i},\tsum_{i=1}^{n}g_{i}(\lambda
_{i}^{g})y_{i})\cap $Max$_{w}\cup _{z\in X}F(x_{i},z)=\emptyset .$

We claim that $F(x_{i},\tsum_{i=1}^{n}g_{i}(\lambda _{i}^{g})y_{i})\cap $Max$%
_{w}\cup _{z\in X}F(x_{i},z)=\emptyset $ for each $i\in \{1,2,...,n\}.$
Indeed, if, by contrary, we assume that there exists $i_{0}\in \{1,2,...,n\}$
and $t\in F(x_{i_{0}},\tsum_{i=1}^{n}g_{i}(\lambda _{i}^{g})y_{i})$ such
that $t\in $Max$_{w}\cup _{z\in X}F(x_{i_{0}},z),$ then, it is true that $%
t\in F(\tsum_{i=1}^{n}\lambda _{i}^{g}x_{i},\tsum_{i=1}^{n}g_{i}(\lambda
_{i}^{g})y_{i})-S$ (1) and $t\in $Max$_{w}\cup _{z\in X}F(x_{i_{0}},z)$ (2)$%
. $

According to (1), we have $t=t^{\prime }-s_{0},$ where $t^{\prime }\in
F(\tsum_{i=1}^{n}\lambda _{i}^{g}x_{i},\tsum_{i=1}^{n}g_{i}(\lambda
_{i}^{g})y_{i})$ and $s_{0}\in S,$ therefore $t^{\prime }=t+s_{0}\in
F(\tsum_{i=1}^{n}\lambda _{i}^{g}x_{i},\tsum_{i=1}^{n}g_{i}(\lambda
_{i}^{g})y_{i}).$ According to the relation (2), $\cup _{z\in
X}F(x_{i_{0}},z)\cap \{t+$int$S\}=\emptyset .$ Consequently, $t^{\prime
}+s\notin \cup _{z\in X}F(x_{i_{0}},z)$ if $s\in $int$S$ (we take into
account that $t^{\prime }+s=t+(s_{0}+s)\in t+$int$S).$ Then, $\cup _{z\in
X}F(x_{i_{0}},z)\cap \{t^{\prime }+$int$S\}=\emptyset ,$ which implies $%
t^{\prime }\in $Max$_{w}\cup _{z\in X}F(x_{i_{0}},z).$

Thus, we have that $t^{\prime }\in F(\tsum_{i=1}^{n}\lambda
_{i}^{g}x_{i},\tsum_{i=1}^{n}g_{i}(\lambda _{i}^{g})y_{i})\cap $Max$_{w}\cup
_{z\in X}F(x_{i_{0}},z),$ which is a contradiction. It remains that $%
F(x_{i},\tsum_{i=1}^{n}g_{i}(\lambda _{i}^{g})y_{i})\cap $Max$_{w}\cup
_{z\in X}F(x_{i},z)=\emptyset $ for each $i\in \{1,2,...,n\}.$

Since $F$ is type-(iii)\textit{\ }pair properly quasi-concave in the second
argument on $X\times X,$ there exists $j\in \{1,2,...,n\}$ such that $%
F(x_{j},y_{j})\cap $Max$_{w}\cup _{z\in X}F(x_{j},z)=\emptyset ,$ which
contradicts the assumption about $(x_{j},y_{j})$. According toTheorem 2.1,
there exists $x^{\ast }\in T(x^{\ast }),$ that is, $F(x^{\ast },x^{\ast
})\cap $Max$_{w}\tbigcup_{y\in X}F(x^{\ast },y)\neq \emptyset .$

(ii) Let us define the set-valued map $Q:X\rightrightarrows X$ by

$Q(y)=\{x\in X:F(x,y)\cap $Min$_{w}\cup _{x\in X}F(x,y)\neq \emptyset \}$
for each $y\in X.$

Further, the proof follows a similar line as above and we conclude that
there exists $y^{\ast }\in Q(y^{\ast }),$ that is, $F(y^{\ast },y^{\ast
})\cap $Min$_{w}\cup _{x\in X}F(x,y^{\ast })\neq \emptyset .$ $\ \ \ \ \ \ \
\ \ \ \ \ \ \ \ \ \ \ \ \ \ \ \ \ \ \ \ \ \ \ \ \ \ \ \ \ \ \ \ \ \ \ \
\square $

\textit{Remark 4.1. }The $S-$transfer type-(v)\textit{\ }$\mu -$convexity of 
$F$ in the first argument on $X\times X$ is verified by all real-valued set
valued maps which fulfill the property that $\tbigcup_{y\in X}F(x,y)$\textit{%
\ }is a compact set for each $x\in X$. This fact is a consequence of Remark
3.1.

As a first application of the previous lemma, we obtain the following
result, which differs from Theorem 3.1 in [32] becose we only take into
consideration the hypothesis which concern convexity properties of
set-valued maps. No form of continuity is assumed.

\textbf{Theorem 4.1} \textit{Let }$X$\textit{\ be a (n-1) dimensional
simplex of a Hausdorff topological vector space }$E,Y$\textit{\ be a compact
set in a Hausdorff topological vector space }$Z$\textit{\ and let }$S$%
\textit{\ be a pointed closed convex cone in }$Z$\textit{\ with its interior
int}$S\neq \emptyset .$\textit{\ Let }$F:X\times X\rightrightarrows Y$%
\textit{\ be a set-valued map with non-empty values. }

\textit{i) Suppose that }$\tbigcup_{y\in X}F(x,y)$\textit{\ is a compact set
for each }$x\in X$. \textit{If the set-valued map }$F$ \textit{is }$S-$%
\textit{transfer type-(v) }$\mu -$\textit{convex in the first argument on} $%
X\times X$, \textit{type-(iii) pair properly quasi-concave in the second
argument on }$X\times X$ \textit{and }$F(\cdot ,y)$\textit{\ is type-(iii)
naturally }$S-$\textit{quasi-concave on }$X$ \textit{for each }$y\in X,$%
\textit{\ then, there exist the elements }$z_{1}\in $Max$\overline{\cup
_{x\in X}F(x,x)}$\textit{\ and }$z_{2}\in $Min$\overline{\cup _{x\in X}\text{%
Max}_{w}F(x,X)}$ \textit{such that }$z_{1}\in z_{2}+S.$

\textit{ii) Suppose that }$\tbigcup_{x\in X}F(x,y)$\textit{\ is a compact
set for each }$y\in X$. \textit{If the set-valued map }$F$ \textit{is }$S-$%
\textit{transfer type-(v) }$\mu -$\textit{concave in the second argument on} 
$X\times X,$\textit{\ type-(iii) pair properly quasi-convex in the first
argument on }$X\times X$ \textit{and }$F(x,\cdot )$\textit{\ is type-(iii)
naturally }$S-$\textit{quasi-convex\ on }$X$ \textit{for each }$x\in X,$\ 
\textit{then, there exist the elements }$z_{1}\in $Min$\overline{\cup _{x\in
X}F(x,x)}$\textit{\ and }$z_{2}\in $Max$\overline{\cup _{y\in X}\text{Min}%
_{w}F(X,y)}$ \textit{such that }$z_{1}\in z_{2}-S.$

\textit{Proof.} i) According to Lemma 4.1, there exists $x^{\ast }\in X$
such that $F(x^{\ast },x^{\ast })\cap $

Max$_{w}\cup _{y\in X}F(x^{\ast },y)\neq \emptyset .$

We have $F(x^{\ast },x^{\ast })\subset \overline{\cup _{x\in X}F(x,x)}$ and,
according to Lemma 2.1, it follows that $\overline{\cup _{x\in X}F(x,x)}%
\subset $Max $\overline{\cup _{x\in X}F(x,x)}-S,$ so that, $F(x^{\ast
},x^{\ast })\subset $Max$\overline{\cup _{x\in X}F(x,x)}-S.$

On the other hand, Max$_{w}\cup _{y\in X}F(x^{\ast },y)\subset \overline{%
\cup _{x\in X}\text{Max}_{w}F(x,X)}$ and, according to Lemma 2.1, it follows
that $\overline{\cup _{x\in X}\text{Max}_{w}F(x,X)}\subset $Min$\overline{%
\cup _{x\in X}\text{Max}_{w}F(x,X)}+S,$ so that, Max$_{w}\cup _{y\in
X}F(x^{\ast },y)\subset $Min$\overline{\cup _{x\in X}\text{Max}_{w}F(x,X)}%
+S. $

Hence, for every $u\in F(x^{\ast },x^{\ast })$ and $v\in $Max$_{w}\cup
_{y\in X}F(x^{\ast },y),$ there exist the elements $z_{1}\in $Max$\overline{%
\cup _{x\in X}F(x,x)}$ and $z_{2}\in $Min$\overline{\cup _{x\in X}\text{Max}%
_{w}F(x,X)}$ such that $u\in z_{1}-S$ and $v\in z_{2}+S.$ If we take $u=v,$
we have $z_{1}\in z_{2}+S.$ \ \ \ \ \ \ \ \ \ \ \ \ \ \ \ \ \ \ \ \ \ \ \ \
\ \ \ \ \ \ \ \ \ \ \ \ \ \ \ \ \ \ \ \ \ \ \ \ \ \ \ \ \ \ \ \ \ \ \ \ \ \ 

ii) According to Lemma 4.1, there exists $y^{\ast }\in X$ such that $%
F(y^{\ast },y^{\ast })\cap $

Min$_{w}\cup _{x\in X}F(x,y^{\ast })\neq \emptyset .$

We have $F(y^{\ast },y^{\ast })\subset \overline{\cup _{x\in X}F(x,x)}$ and,
according to Lemma 2.1, it follows that $\overline{\cup _{x\in X}F(x,x)}%
\subset $Min $\overline{\cup _{x\in X}F(x,x)}+S,$ so that, $F(y^{\ast
},y^{\ast })\subset $Min$\overline{\cup _{x\in X}F(x,x)}+S.$

On the other hand, Min$_{w}\cup _{x\in X}F(x,y^{\ast })\subset \overline{%
\cup _{y\in X}\text{Min}_{w}F(X,y)}$ and, according to Lemma 2.1, it follows
that $\overline{\cup _{y\in X}\text{Min}_{w}F(X,y)}\subset $Max$\overline{%
\cup _{y\in X}\text{Min}_{w}F(X,y)}-S,$ consequently, Min$_{w}\cup _{x\in
X}F(x,y^{\ast })\subset $Max$\overline{\cup _{y\in X}\text{Min}_{w}F(X,y)}%
-S. $

Hence, for every $u\in F(y^{\ast },y^{\ast })$ and $v\in $Min$_{w}\cup
_{x\in X}F(x,y^{\ast }),$ there exist the elements $z_{1}\in $Min$\overline{%
\cup _{x\in X}F(x,x)}$ and $z_{2}\in $Max$\overline{\cup _{y\in X}\text{Min}%
_{w}F(X,y)}$ such that $u\in z_{1}+S$ and $v\in z_{2}-S.$ If we take $u=v,$
we have $z_{1}\in z_{2}-S.$ \ \ \ \ \ \ \ \ \ \ \ \ \ \ \ \ \ \ \ \ \ \ \ \
\ \ \ \ \ \ \ \ \ \ \ \ \ \ \ \ \ \ \ \ $\square $\medskip

An important version of Theorem 4.1 is obtained in the case when the
set-valued map has the property $\alpha $ (resp.$\alpha ^{\prime }).$

\textbf{Theorem 4.2} \textit{Let }$X$\textit{\ be a (n-1) dimensional
simplex of a Hausdorff topological vector space }$E,Y$\textit{\ be a compact
set in a Hausdorff topological vector space }$Z$\textit{\ and let }$S$%
\textit{\ be a pointed closed convex cone in }$Z$\textit{\ with its interior
int}$S\neq \emptyset .$\textit{\ Let }$F:X\times X\rightrightarrows Y$%
\textit{\ be a set-valued map with nonempty values. }

\textit{i) Suppose that }$F$ \textit{satisfies the property} $\alpha .$ 
\textit{If }$F$ \textit{is type-(iii) pair properly quasi-concave in the
second argument on }$X\times X$ \textit{and }$F(\cdot ,y)$\textit{\ is
type-(iii) naturally }$S-$\textit{quasi-concave on }$X$ \textit{for each }$%
y\in X,$\textit{\ then, there exist the elements }$z_{1}\in $Max$\overline{%
\cup _{x\in X}F(x,x)}$\textit{\ and }$z_{2}\in $Min$\overline{\cup _{x\in X}%
\text{Max}_{w}F(x,X)}$ \textit{such that }$z_{1}\in z_{2}+S.$

\textit{ii) Suppose that }$F$ \textit{satisfies the property} $\alpha
^{\prime }.$ \textit{If }$F$ \textit{is\ type-(iii) pair properly
quasi-convex in the first argument on }$X\times X$ \textit{and }$F(x,\cdot )$%
\textit{\ is type-(iii) naturally }$S-$\textit{quasi-convex\ on }$X$ \textit{%
for each }$x\in X,$\ \textit{then, there exist the elements }$z_{1}\in $Min$%
\overline{\cup _{x\in X}F(x,x)}$\textit{\ and }$z_{2}\in $Max$\overline{\cup
_{y\in X}\text{Min}_{w}F(X,y)}$ \textit{such that }$z_{1}\in z_{2}-S.$

\textit{Example 4.1} Let $S=-R_{+}^{2}$, and for each $x\in \lbrack 0,1],$
let $S^{\ast }((0,0),x)=\{(u,v)\in \lbrack 0,1]\times \lbrack
0,1]:u^{2}\times v^{2}\leq x^{2}\}$ and $F:[0,1]\times \lbrack
0,1]\rightrightarrows \lbrack 0,1]\times \lbrack 0,1]$ be defined by

$F(x,y)=\left\{ 
\begin{array}{c}
\{(0,0)\}\text{ \ \ \ \ for each\ \ \ }0\leq x\leq y\leq 1; \\ 
S^{\ast }((0,0),x)\text{ for each }0\leq y<x\leq 1.%
\end{array}%
\right. $

We notice that $F$ is not continuous on $X.$

a) $F$ is $-R_{+}^{2}-$transfer type-(v)\textit{\ }$\mu -$convex in the
first argument$.$

Let $x_{1},x_{2},...,x_{n}\in \lbrack 0,1]$ and $z\in \lbrack 0,1].$ For
each $i\in \{1,2,...,n\},$ there exists $%
z_{i}=z_{i}(x_{1},x_{2},...,x_{n},z)\geq \max_{i=1,2,...,n}x_{i}\in \lbrack
0,1]$ such that $F(x_{i},z_{i})=\{(0,0)\}$ for each $i=1,2,...n$ and then, $%
F(\tsum_{i=1}^{n}\lambda _{i}x_{i},z)\cap (\tbigcup_{y\in
X}F(x_{i},y))\subset \{(0,0)\}-(-R_{+}^{2})$ for each $\lambda =(\lambda
_{1},\lambda _{2},...,\lambda _{n})\in \Delta _{n-1}$.

It follows that $F$ is $-R_{+}^{2}-$transfer type-(v) $\mu -$convex in the
first argument on $[0,1]\times \lbrack 0,1].$

b) $F$ is type-(iii)\textit{\ }pair properly $-R_{+}^{2}-$quasiconcave in
the second argument on $[0,1]\times \lbrack 0,1].$

Let us consider $(x_{1},y_{1}),$ $(x_{2},y_{2})\in \lbrack 0,1]\times
\lbrack 0,1]$ and let us assume, without loss of generalization, that $%
y_{1}\leq y(\lambda )\leq y_{2}$ for each $\lambda \in \lbrack 0,1],$ where $%
y(\lambda )=\lambda y_{1}+(1-\lambda )y_{2}.$

$F(x_{1},y_{1})=\left\{ 
\begin{array}{c}
\{(0,0)\}\text{ \ \ \ \ \ \ \ \ \ \ \ \ \ \ for each\ \ \ \ \ \ \ \ \ \ \ }%
0\leq x_{1}\leq y_{1}\leq 1; \\ 
S^{\ast }((0,0),x_{1})\text{ \ \ \ \ \ \ \ \ for each \ \ \ \ \ \ \ \ }0\leq
y_{1}<x_{1}\leq 1,%
\end{array}%
\right. $

$F(x_{2},y_{2})=\left\{ 
\begin{array}{c}
\{(0,0)\}\text{ \ \ \ for each\ \ \ \ \ \ }0\leq x_{2}\leq y_{2}\leq 1; \\ 
S^{\ast }((0,0),x_{2})\text{ \ for each }0\leq y_{2}<x_{2}\leq 1%
\end{array}%
\right. $,

$F(x_{1},y(\lambda ))=\left\{ 
\begin{array}{c}
\{(0,0)\}\text{\ \ \ \ \ for each\ \ \ \ \ }0\leq x_{1}\leq y(\lambda )\leq
1; \\ 
S^{\ast }((0,0),x_{1})\text{ for each }0\leq y(\lambda )<x_{1}\leq 1;%
\end{array}%
\right. $

$F(x_{2},y(\lambda ))=\left\{ 
\begin{array}{c}
\{(0,0)\}\text{ \ \ \ \ for each\ \ \ \ }0\leq x_{2}\leq y(\lambda )\leq 1;
\\ 
S^{\ast }((0,0),x_{2})\text{ for each }0\leq y(\lambda )<x_{2}\leq 1.%
\end{array}%
\right. $

b1) If $x_{1}\leq y_{1}\leq y(\lambda ),$ then, $F(x_{1},y_{1})=\{(0,0)\},$ $%
F(x_{1},y(\lambda ))=\{(0,0)\};$

b2) if $y_{1}\leq y(\lambda )<x_{1},$ then, $F(x_{1},y_{1})=S^{\ast
}((0,0),x_{1}),$ $F(x_{1},y(\lambda ))=S^{\ast }((0,0),x_{1});$

b3) if $y_{1}<x_{1}\leq y(\lambda ),$ then, $F(x_{1},y_{1})=S^{\ast
}((0,0),x_{1}),$ $F(x_{1},y(\lambda ))=\{(0,0)\}.$

Then, $F(x_{1},y_{1})\subset F(x_{1},y(\lambda ))-(-R_{+}^{2})$ for each $%
\lambda \in \lbrack 0,1].$

c) We prove that\ $F(\cdot ,y)$\ is type-(iii)\textit{\ }naturally $%
-R_{+}^{2}-$quasiconcave on $[0,1]$ for each $y\in \lbrack 0,1].$

Let $y\in \lbrack 0,1]$ be fixed, $x_{1},x_{2}\in \lbrack 0,1]$, $\lambda
\in \lbrack 0,1]$ and $x(\lambda )=\lambda x_{1}+(1-\lambda )x_{2}$.

c1) If $x_{1}\geq x_{2}>y,$ $F(x_{1},y)=S^{\ast }((0,0),x_{1}),$ $%
F(x_{2},y)=S^{\ast }((0,0),x_{2})$, $F(x(\lambda ),y)=S^{\ast
}((0,0),x(\lambda ))$ and

co$\{F(x_{1},y)\cup F(x_{2},y)\}=S^{\ast }((0,0),x_{1})\subset S^{\ast
}((0,0),x(\lambda ))-(-R_{+}^{2})=F(x(\lambda ),y)-(-R_{+}^{2});$

c2) if $x_{1}\leq x_{2}\leq y,$ $F(x_{1},y)=\{(0,0)\},$ $F(x_{2},y)=\{(0,0)%
\} $, $F(x(\lambda ),y)=\{(0,0)\}$ and co$\{F(x_{1},y)\cup
F(x_{2},y)\}=\{(0,0)\}\subset F(x(\lambda ),y)-(-R_{+}^{2});$

c3) if $x_{1}>y\geq x_{2},$ then, $F(x_{1},y)=S^{\ast }((0,0),x_{1}),$ $%
F(x_{2},y)=\{(0,0)\}$ and

co$\{F(x_{1},y)\cup F(x_{2},y)\}=S^{\ast }((0,0),x_{1});$

\ \ \ \ if $x_{1}\geq x(\lambda )>y\geq x_{2},$ then, $F(x(\lambda
),y)=S^{\ast }((0,0),x(\lambda ))$ and

co$\{F(x_{1},y)\cup F(x_{2},y)\}=S^{\ast }((0,0),x_{1})\subset F(x(\lambda
),y)-(-R_{+}^{2});$

\ \ \ \ \ if $x_{1}>y\geq x(\lambda )\geq x_{2},$ then, $F(x(\lambda
),y)=\{(0,0)\}$ and

co$\{F(x_{1},y)\cup F(x_{2},y)\}=S^{\ast }((0,0),x_{1})\subset
\{(0,0)\}-(-R_{+}^{2})=F(x(\lambda ),y)-(-R_{+}^{2}).$

The following equality is true:

$\cup _{y\in X}F(x,y)=S^{\ast }((0,0),x)$ and, consequently, $\cup _{y\in
X}F(x,y)$ is a compact set.

All the assumptions of Theorem 4.2 are fulfilled, then, there exist the
elements $z_{1}\in $Max$\overline{\cup _{x\in X}F(x,x)}$\ and $z_{2}\in $Min$%
\overline{\cup _{x\in X}\text{Max}_{w}F(x,X)}$ such that $z_{1}\in
z_{2}+(-R_{+}^{2}).$

In our case, $\cup _{x\in X}F(x,x)=\{(0,0)\},$ Max$\overline{\cup _{x\in
X}F(x,x)}=\{(0,0)\},$ Max$_{w}F(x,X)=\{(0,0)\}$ and Min$\overline{\cup
_{x\in X}\text{Max}_{w}F(x,X)}=\{(0,0)\}.$ Then, taking $z_{1}=(0,0)$ and $%
z_{2}=(0,0),$ we have that $z_{1}\in z_{2}+(-R_{+}^{2}).$ \ \ \ \medskip

Considering Remark 4.2, we obtain the following result as a consequence of
Theorem 4.2, for the real-valued maps case.

\textbf{Corollary 4.1} \textit{Let }$X$\textit{\ be a (n-1) dimensional
simplex of a Hausdorff topological vector space }$E,Y$\textit{\ a compact
set in }$\mathit{R}$\textit{\ and let }$S$\textit{\ be a pointed closed
convex cone in }$R$\textit{\ with its interior int}$S\neq \emptyset .$%
\textit{\ Let }$F:X\times X\rightrightarrows Y$\textit{\ be a set-valued map
with non-empty values. }

\textit{i) Suppose that }$\tbigcup_{y\in X}F(x,y)$\textit{\ is a compact set
for each }$x\in X$. \textit{If }$F$\textit{\ is type-(iii) pair properly
quasi-concave in the second argument on }$X\times X$ \textit{and }$F(\cdot
,y)$\textit{\ is type-(iii) naturally }$S-$\textit{quasi-concave on }$X$ 
\textit{for each }$y\in X,$\textit{\ then, there exist }$z_{1}\in $Max$%
\overline{\cup _{x\in X}F(x,x)}$\textit{\ and }$z_{2}\in $Min$\overline{\cup
_{x\in X}\text{Max}_{w}F(x,X)}$ \textit{such that }$z_{1}\in z_{2}+S.$

\textit{ii) Suppose that }$\tbigcup_{x\in X}F(x,y)$\textit{\ is a compact
set for each }$y\in X$. \textit{If} $F$ \textit{is} \textit{type-(iii) pair
properly quasi-convex in the first argument on }$X\times X$ \textit{and }$%
F(x,\cdot )$\textit{\ is type-(iii) naturally }$S-$\textit{quasi-convex\ on }%
$X$ \textit{for each }$x\in X,$\ \textit{then, there exist }$z_{1}\in $Min$%
\overline{\cup _{x\in X}F(x,x)}$\textit{\ and }$z_{2}\in $Max$\overline{\cup
_{y\in X}\text{Min}_{w}F(X,y)}$ \textit{such that }$z_{1}\in z_{2}-S.$%
\medskip

\textit{Example 4.2} Let $X=[0,1],$ $Y=[-1,1],$ $S=[0,\infty )$ and $%
F:X\times X\rightrightarrows Y$ be defined by $F(x,y)=\left\{ 
\begin{array}{c}
\lbrack -1,1]\text{ if }0\leq x\leq y\leq 1; \\ 
\lbrack -x,1]\text{ if }0\leq y<x\leq 1.%
\end{array}%
\right. $

We notice that $F$ is not continuous on $X$ and it is $S-$transfer type-(v)%
\textit{\ }$\mu -$convex in the first argument$.$

a) In Example 3.5 we have seen that $F$ is type-(iii)\textit{\ }pair
properly quasiconcave in the second argument on $X\times X.$

b) We prove that\ $F(\cdot ,y)$\ is type-(iii)\textit{\ }naturally $S-$%
quasiconcave on $X$ for each $y\in X.$

Let $y\in \lbrack 0,1]$ be fixed, $x_{1},x_{2}\in \lbrack 0,1]$, $\lambda
\in \lbrack 0,1]$ and $x(\lambda )=\lambda x_{1}+(1-\lambda )x_{2}$.

b1) If $x_{1}\geq x_{2}\geq y,$ $F(x_{1},y)=[-x_{1},1],$ $%
F(x_{2},y)=[-x_{2},1]$, $F(x(\lambda ),y)=[-x(\lambda ),1]$

and co$\{F(x_{1},y)\cup F(x_{2},y)\}=[-x_{1},1]\subset \lbrack -x(\lambda
),1]-[0,\infty )=F(x(\lambda ),y)-[0,\infty );$

b2) if $x_{1}\leq x_{2}\leq y,$ $F(x_{1},y)=[-1,1],$ $F(x_{2},y)=[-1,1]$, $%
F(x(\lambda ),y)=[-1,1]$ and

co$\{F(x_{1},y)\cup F(x_{2},y)\}=[-1,1]\subset \lbrack -1,1]-[0,\infty
)=F(x(\lambda ),y)-[0,\infty );$

b3) if $x_{1}\geq y\geq x_{2},$ then, $F(x_{1},y)=[-x_{1},1],$ $%
F(x_{2},y)=[-1,1]$ and

co$\{F(x_{1},y)\cup F(x_{2},y)\}=[-1,1];$

\ if $x_{1}\geq x(\lambda )\geq y\geq x_{2},$ then, $F(x(\lambda
),y)=[-x(\lambda ),1]$ and

co$\{F(x_{1},y)\cup F(x_{2},y)\}=[-1,1]\subset \lbrack -x(\lambda
),1]-[0,\infty )=F(x(\lambda ),y)-[0,\infty );$

if $x_{1}\geq y\geq x(\lambda )\geq x_{2},$ then, $F(x(\lambda ),y)=[-1,1]$
and

co$\{F(x_{1},y)\cup F(x_{2},y)\}=[-1,1]\subset \lbrack -1,1]-[0,\infty
)=F(x(\lambda ),y)-[0,\infty ).$

The following equalities are true:

$\cup _{y\in X}F(x,y)=\cup _{y<x}[-x,1]\cup \cup _{y\geq x}[-1,1]=[-x,1]\cup
\lbrack -1,1]=[-1,1]$ and, consequently, $\cup _{y\in X}F(x,y)$ is a compact
set.

All the assumptions of Corollary 4.1 are fulfilled, then, there exist the
elements $z_{1}\in $Max$\overline{\cup _{x\in X}F(x,x)}$\ and $z_{2}\in $Min$%
\overline{\cup _{x\in X}\text{Max}_{w}F(x,X)}$ such that $z_{1}\in z_{2}+S.$

In our case, $\cup _{x\in X}F(x,x)=[-1,1],$ Max$\overline{\cup _{x\in
X}F(x,x)}=\{1\},$ Max$_{w}F(x,X)=\{1\}$ and Min$\overline{\cup _{x\in X}%
\text{Max}_{w}F(x,X)}=\{1\}.$ Then, taking $z_{1}=1$ and $z_{2}=1,$ we have
that $z_{1}\in z_{2}+S.$ \ \ \ \medskip

The next corollary is a particular case of Theorem 4.1.

\textbf{Corollary 4.2} \textit{Let }$X$\textit{\ be a (n-1) dimensional
simplex of a Hausdorff topological vector space }$E,$\textit{\ }$Y$\textit{\
be a compact set in a Hausdorff topological vector space }$Z$\textit{\ and
let }$S$\textit{\ be a pointed closed convex cone in }$Z$\textit{\ with its
interior int}$S\neq \emptyset .$\textit{\ Let }$f:X\times X\rightarrow Y$%
\textit{\ be a vector-valued mapping.}

\textit{i)} \textit{Suppose that }$\tbigcup_{y\in X}f(x,y)$\textit{\ is a
compact set for each }$x\in X$. \textit{If the mapping }$f$ \textit{is }$S-$%
\textit{transfer }$\mu -$\textit{convex in the first argument on} $X\times X$%
,\textit{\ pair properly quasi-concave in the second argument on }$X\mathit{%
\times X}$ \textit{and }$f(\cdot ,y)$\textit{\ is naturally }$S-$\textit{%
quasi-concave on }$X$ \textit{for each }$y\in X,$\textit{\ then, there exist 
}$z_{1}\in $Max$\overline{\cup _{x\in X}f(x,x)}$\textit{\ and }$z_{2}\in $Min%
$\overline{\cup _{x\in X}\text{Max}_{w}f(x,X)}$ \textit{such that }$z_{1}\in
z_{2}+S.$

\textit{ii)} \textit{Suppose that }$\tbigcup_{x\in X}f(x,y)$\textit{\ is a
compact set for each }$y\in X$. \textit{If the mapping }$f$ \textit{is }$S$ 
\textit{transfer }$\mu -$\textit{concave in the second argument on} $X\times
X$\textit{, pair properly quasi-convex in the first argument on }$X\times X$ 
\textit{and \ }$f(x,\cdot )$\textit{\ is naturally }$S$ \textit{quasi-convex
on }$X$ \textit{for each }$x\in X,$\textit{\ then, there exist }$z_{1}\in $%
Min$\overline{\cup _{x\in X}f(x,x)}$\textit{\ and }$z_{2}\in $Max$\overline{%
\cup _{y\in X}\text{Min}_{w}f(X,y)}$\textit{such that }$z_{1}\in z_{2}-S.$%
\medskip

We search to weaken the assumptions from Lemma 4.1, especially the $S-$%
transfer $\mu -$convexity (resp. $S-$transfer $\mu -$concavity) and the
naturally $S-$quasi-concavity (resp. naturally $S-$quasi-convexity), but
another proving method needs to be used: we build a constant selection for a
set-valued map. This change requires a new condition instead of pair
quasi-convexity (resp. pair quasi-concavity), a condition we called $\gamma $
(resp. $\gamma ^{\prime }$). Under the condition $\gamma $ (resp.$\gamma
^{\prime }$)$,$ the assumption of transfer properly $S-$quasi-concavity
(resp. transfer properly $S-$quasi-convexity) proves to be necessary. The
next Lemma is the key used in order to obtain Theorem 4.3.

\textbf{Lemma 4.2} \textit{Let }$X$\textit{\ be a convex set in a Hausdorff
topological vector space }$E,$\textit{\ }$Y$\textit{\ a compact set in the
Hausdorff topological vector space }$Z$\textit{\ and let }$S$\textit{\ a
pointed closed convex cone in }$Z$\textit{\ with its interior int}$S\neq
\emptyset .$\textit{\ Let }$F:X\times X\rightrightarrows Y$\textit{\ be a
set-valued map with non-empty values.}

\textit{(i) Suppose that }$\tbigcup_{y\in X}F(x,y)$\textit{\ is a compact
set for each }$x\in X$. \textit{If }$F$ \textit{is }$S-$\textit{transfer
weakly type-(v) }$\mu -$\textit{convex in the first argument on} $X\times X$,%
\textit{\ transfer type-(iii) properly }$S-$\textit{quasi-concave in the
first argument on }$X\times X$ \textit{and satisfies the condition} $\gamma
, $\textit{\ then there exists }$x^{\ast }\in X$\textit{\ such that }$%
F(x^{\ast },x^{\ast })\cap $Max$_{w}\tbigcup_{y\in X}F(x^{\ast },y)\neq
\emptyset .$

\textit{(ii) Suppose that }$\tbigcup_{x\in X}F(x,y)$\textit{\ is a compact
set for each }$y\in X$. \textit{\ If }$F$ \textit{is transfer weakly
type-(v) }$\mu -$\textit{concave in the second argument on }$X\times X$%
\textit{, transfer type-(iii) properly }$S-$\textit{quasi-convex\ in the
second argument on }$X\times X,$\textit{\ and satisfies the condition} $%
\gamma ^{\prime }$\textit{, then, there exists }$y^{\ast }\in X$\textit{\
such that }$F(y^{\ast },y^{\ast })\cap $Min$_{w}\cup _{x\in X}F(x,y^{\ast
})\neq \emptyset .\medskip $

\textit{Proof.} Let us define the set-valued map $T:X\rightrightarrows X$ by

$T(x)=\{y\in X:F(x,y)\cap $Max$_{w}\cup _{z\in X}F(x,z)\neq \emptyset \}$
for each $x\in X.$

We claim that $T$ is non-empty valued. Indeed, since $\cup _{z\in X}F(x,z)$
is a compact set for each $x\in X,$ by Lemma 2.1, Max$_{w}\cup _{z\in
X}F(x,z)\neq \emptyset .$ For each $x\in X,$ let $z_{x}\in $Max$_{w}\cup
_{z\in X}F(x,z).$ Then, there exists $y_{x}\in X$ such that $z_{x}\in
F(x,y_{x}).$ It is clear that, $y_{x}\in T(x)=\{y\in X:F(x,y)\cap $Max$%
_{w}\cup _{z\in X}F(x,z)\neq \emptyset \}$ and consequently, $T(x)\neq
\emptyset $ for each $x\in X.$

Since $F$ satisfies the condition $\gamma $, there exist $n\in N,$ $%
(x_{1},y_{1}),(x_{2},y_{2})...,(x_{n},y_{n})\in X\times X$ and $y^{\ast }\in 
$co$\{x_{i}:i=1,2,...,n\}$ such that $F(x_{i},y_{i})\subset F(x_{i},y^{\ast
})-S$ and $F(x_{i},y_{i})\cap $Max$_{w}\cup _{z\in X}F(x_{i},z)\neq
\emptyset $ for each $i\in \{1,2,...,n\}.$

Let us fix $i_{0}\in \{1,2,...,n\}.$ There exists $t_{i_{0}}\in
F(x_{i_{0}},y_{i_{0}})\cap $Max$_{w}\cup _{z\in X}F(x_{i_{0}},z).$ This
means that $t_{i_{0}}\in F(x_{i_{0}},y_{i_{0}})$ and $\cup _{z\in
X}F(x_{i_{0}},z)\cap (t_{i_{0}}+$int$S)=\emptyset .$ There exists $%
t_{i_{0}}^{\prime }\in F(x_{i_{0}},y^{\ast })$ and $s_{i_{0}}\in S$ such
that $t_{i_{0}}^{\prime }=t_{i_{0}}+s_{i_{0}}.$ Therefore, $%
t_{i_{0}}^{\prime }\in \cup _{z\in X}F(x_{i_{0}},z)$ and, for each $%
s^{\prime }\in $int$S,$ $(t_{i_{0}}^{\prime }+s^{\prime })\cap \cup _{z\in
X}F(x_{i_{0}},z)=(t_{i_{0}}+s_{i_{0}}+s^{\prime })\cap \cup _{z\in
X}F(x_{i_{0}},z)=\emptyset .$ It follows that $(t_{i_{0}}^{\prime }+$int$%
S)\cap \cup _{z\in X}F(x_{i_{0}},z)=\emptyset $, and, since $%
t_{i_{0}}^{\prime }\in \cup _{z\in X}F(x_{i_{0}},z),$ we have that $%
t_{i_{0}}\in F(x_{i_{0}},y^{\ast })\cap $Max$_{w}\cup _{z\in
X}F(x_{i_{0}},z).$ We showed that $y^{\ast }\in T(x_{i_{0}}),$ and, since $%
i_{0}$ is arbitrary and $y^{\ast }\in $co$\{x_{i}:i=1,2,...,n\}$, then, $%
y^{\ast }\in \tbigcap\limits_{i=1}^{n}T(x_{i})\cap $co$\{x_{i}:i=1,2,...,n%
\}. $ Hence, $\tbigcap\limits_{i=1}^{n}T(x_{i})$ is non-empty.

Further, we will prove that $T$ is quasi-convex. By contrary, we assume that 
$T$ is not quasi-convex. Then, suppose that there exists $z^{\ast }\in
\tbigcap\limits_{i=1}^{n}T(x_{i})$ and $\lambda ^{\ast }\in \Delta _{n-1}$
such that $z^{\ast }\notin T(\tsum_{i=1}^{n}\lambda _{i}^{\ast }x_{i}),$
that is $F(x_{i},z^{\ast })\cap $Max$_{w}\cup _{z\in X}F(x_{i},z)\neq
\emptyset $ for each $i\in \{1,2,...,n\}$ and $F(\tsum_{i=1}^{n}\lambda
_{i}^{\ast }x_{i},z^{\ast })\cap $Max$_{w}\cup _{z\in
X}F(\tsum_{i=1}^{n}\lambda _{i}^{\ast }x_{i},z)=\emptyset .$

Since $F$ is $S-$transfer weakly type-(v)\textit{\ }$\mu -$convex in the
first argument and we also have $F(\tsum_{i=1}^{n}\lambda _{i}^{\ast
}x_{i},z^{\ast })\cap $Max$_{w}\tbigcup_{z\in X}F(\tsum_{i=1}^{n}\lambda
_{i}^{\ast }x_{i},z)=\emptyset $, it follows that, there exists $i_{0}\in I$
and $z_{i_{0}}\in X$ such that $F(\tsum_{i=1}^{n}\lambda _{i}^{\ast
}x_{i},z^{\ast })\cap (\tbigcup_{z\in X}F(x_{i_{0}},z))\subset
F(x_{i_{0}},z_{i_{0}})-$int$S.$

Let $t\in F(\tsum_{i=1}^{n}\lambda _{i}^{\ast }x_{i},z^{\ast })\cap
(\tbigcup_{z\in X}F(x_{i_{0}},z))$ and let $u_{i_{0}}\in
F(x_{i_{0}},z_{i_{o}})$ such that $t=u_{i_{0}}-s_{i_{0}},$ $s_{i_{0}}\in $int%
$S.$ Since $t\in F(x_{i_{0}},z_{i_{0}})-$int$S,$ it follows that $%
u_{i_{0}}\in \tbigcup_{z\in X}F(x_{i_{0}},z)\cap \{t+$int$S\}\neq \emptyset
, $ that is $t\in F(\tsum_{i=1}^{n}\lambda _{i}^{\ast }x_{i},z^{\ast })\cap
(\tbigcup_{z\in X}F(x_{i_{0}},z))$ implies the fact that $t\notin $Max$%
_{w}\cup _{z\in X}F(x_{i_{0}},z).$

Consequently, $F(\tsum_{i=1}^{n}\lambda _{i}^{\ast }x_{i},z^{\ast })\cap $Max%
$_{w}\cup _{z\in X}F(x_{i_{0}},z)=\emptyset .$

We claim that $F(x_{i_{0}},z^{\ast })\cap $Max$_{w}\cup _{z\in
X}F(x_{i_{0}},z)=\emptyset .$ On the contrary, we assumethat there exists $%
t\in F(x_{i_{0}},z^{\ast })$ such that $t\in $Max$_{w}\cup _{z\in
X}F(x_{i_{0}},z).$ Since $F$ is transfer type-(iii)\textit{\ }properly $S-$%
quasi-concave in the first argument, then, $t\in F(\tsum_{i=1}^{n}\lambda
_{i}^{\ast }x_{i},z^{\ast })-S.$ We have $t=t^{\prime }-s_{0},$ where $%
t^{\prime }\in F(\tsum_{i=1}^{n}\lambda _{i}^{\ast }x_{i},z^{\ast })$ and $%
s_{0}\in S,$ therefore $t^{\prime }=t+s_{0}\in F(\tsum_{i=1}^{n}\lambda
_{i}^{\ast }x_{i},z^{\ast })$. Since $t\in $Max$_{w}\cup _{z\in
X}F(x_{i_{0}},z),$ $F(x_{i_{0}},z)\cap \{t+$int$S\}=\emptyset .$ For each $%
s\in $int$S,$ $t^{\prime }+s=t+s_{0}+s\in t+$int$S$, which implies $%
t^{\prime }+s\notin F(x_{i_{0}},z)$, that is, $F(x_{i_{0}},z)\cap
\{t^{\prime }+$int$S\}=\emptyset ,$ or, equivalently, $t^{\prime }\in $Max$%
_{w}\cup _{z\in X}F(x_{i_{0}},z).$ We obtained $t^{\prime }\in
F(\tsum_{i=1}^{n}\lambda _{i}^{\ast }x_{i},z^{\ast })\cap $Max$_{w}\cup
_{z\in X}F(x_{i_{0}},z),$ which is a contradiction. It remains that $%
F(x_{i_{0}},z^{\ast })\cap $Max$_{w}\cup _{z\in X}F(x_{i_{0}},z)=\emptyset ,$
and then, $z^{\ast }\notin T(x_{i_{0}}),$ which contradicts $z^{\ast }\in
\tbigcap\limits_{i=1}^{n}T(x_{i})$. Therefore, $T$ is quasi-convex.

We proved that there exist the elements $x^{\ast },x_{1},x_{2},...,x_{n}\in
X $ such that $x^{\ast }\in \tbigcap\limits_{i=1}^{n}T(x_{i})\cap $co$%
\{x_{i}:i=1,2,...,n\}\subset T(x)$ for each $x\in $co$\{x_{i}:i=1,2,..,n\}$,
then, $x^{\ast }\in T(x^{\ast }),$ that is, $F(x^{\ast },x^{\ast })\cap $Max$%
_{w}\tbigcup_{y\in X}F(x^{\ast },y)\neq \emptyset .$

(ii) Let us define the set-valued map $Q:X\rightrightarrows X$ by

$Q(y)=\{x\in X:F(x,y)\cap $Min$_{w}\cup _{x\in X}F(x,y)\neq \emptyset \}$
for each $y\in X.$

Further, the proof follows a similar line as above and we conclude that
there exists $y^{\ast }\in Q(y^{\ast }),$ that is, $F(y^{\ast },y^{\ast
})\cap $Min$_{w}\tbigcup_{x\in X}F(x,y^{\ast })\neq \emptyset .\ \ \ \ \ \ \
\ \ \ \ \ \ \ \ \ \ \ \ \ \ \ \ \ \ \ \ \ \ \ \ \ \ \ \ \ \ \ \ \ \ \ \
\square \medskip $

\textit{Remark 4.2. }The\textit{\ }transfer type-(iii)\textit{\ }properly $%
S- $quasiconcavity in the first argument of $F$ is a necessary condition for
Lemma 4.2 i). In the following example, we have that $F$ satisfies the
condition $\gamma $, it is not transfer type-(iii)\textit{\ }properly $S-$%
quasiconcave in the first argument and the conclusion of Lemma 4.2 i) is not
fulfilled.

Let $X=[0,1],$ $Y=[0,1]$ and $F:X\times X\rightrightarrows Y$ be defined by

$F(x,y)=\left\{ 
\begin{array}{c}
\lbrack 0,1]\text{ if }(x,y)\in \lbrack \frac{1}{4},\frac{3}{4}]\times
\{1\}\cup ([0,\frac{1}{4}]\cup \lbrack \frac{3}{4},1])\times \{\frac{1}{2}\};
\\ 
\{(0\}\text{ \ \ \ \ \ \ \ \ \ \ \ \ \ \ \ \ \ \ \ \ \ \ \ \ \ \ \ \ \ \ \ \
\ \ \ \ \ \ \ \ \ \ \ \ \ \ \ otherwise.}%
\end{array}%
\right. \medskip $

\textit{Remark 4.3. }The two assumptions from Lemma 4.2 i), namely, the $S-$%
transfer weakly type-(v)\textit{\ }$\mu -$convexity in the first argument
and the transfer type-(iii)\textit{\ }properly $S-$quasiconcavity of $F$ in
the first argument on $X\times X,$ imply the following:

for each $x_{1},x_{2},...,x_{n}\in X$ and $z\in X$ $,$ there exists $\lambda
^{\ast }\in \Delta _{n-1},$ $i_{o}\in \{1,2,...n\}$ and $z_{i_{0}}\in X$
such that if $F(\tsum_{i=1}^{n}\lambda _{i}^{\ast }x_{i},z)\cap $Max$%
_{y}\tbigcup_{y\in X}F(x_{i_{0}},y)=\emptyset ,$ it follows that

$F(\tsum_{i=1}^{n}\lambda _{i}^{\ast }x_{i},z)\cap \tbigcup\limits_{y\in
X}F(x_{i_{0}},y)\subset F(\tsum_{i=1}^{n}\lambda _{i}^{\ast
}x_{i},z_{i_{0}}).\medskip $

As a first application of the previous lemma, we obtain the following
result, which differs from Theorem 3.1 in [32] by the fact that the
continuity assumptions are dropped.

\textbf{Theorem 4.3} \textit{Let }$X$ \textit{be a convex set}\ \textit{be
in a Hausdorff topological vector space }$E,$ $Y$\textit{\ a compact set in
a Hausdorff topological vector space }$Z$\textit{\ and let }$S$\textit{\ be
a pointed closed convex cone in }$Z$\textit{\ with its interior int}$S\neq
\emptyset .$\textit{\ Let }$F:X\times X\rightrightarrows Y$\textit{\ be a
set-valued map with non-empty values.}

\textit{i) Suppose that }$\tbigcup_{y\in X}F(x,y)$\textit{\ is a compact set
for each }$x\in X$. \textit{If }$F$ \textit{is }$S-$\textit{transfer weakly
type-(v) }$\mu -$\textit{convex in the first argument on} $X\times X$,%
\textit{\ transfer type-(iii) properly }$S-$\textit{quasi-concave in the
first argument on }$X\times X$ \textit{and} \textit{satisfies the condition} 
$\gamma ,$\textit{\ then, there exist }$z_{1}\in $Max$\overline{\cup _{x\in
X}F(x,x)}$\textit{\ and }$z_{2}\in $Min$\overline{\cup _{x\in X}\text{Max}%
_{w}F(x,X)}$ \textit{such that }$z_{1}\in z_{2}+S.$

\textit{ii)} \textit{Suppose that }$\tbigcup_{x\in X}F(x,y)$\textit{\ is a
compact set for each }$y\in X$. \textit{If }$F$ \textit{is transfer weakly
type-(v) }$\mu -$\textit{concave in the second argument on} $X\times X$,%
\textit{\ transfer type-(iii) properly }$S-$\textit{quasi-convex\ in the
second argument on }$X\times X$ \textit{and} \textit{satisfies the condition}
$\gamma ^{\prime },$\textit{\ then there exist }$z_{1}\in $Min$\overline{%
\cup _{x\in X}F(x,x)}$\textit{\ and }$z_{2}\in $Max$\overline{\cup _{y\in X}%
\text{Min}_{w}F(X,y)\text{ }}$\textit{such that }$z_{1}\in z_{2}-S.$\medskip

\textit{Proof.} i) According to Lemma 4.2, in the case i) there exists $%
x^{\ast }\in X$ such that $F(x^{\ast },x^{\ast })\cap $Max$_{w}\cup _{y\in
X}F(x^{\ast },y)\neq \emptyset $ and in the case ii), there exists $y^{\ast
}\in X$ such that $F(y^{\ast },y^{\ast })\cap $Min$_{w}\cup _{x\in
X}F(x,y^{\ast })\neq \emptyset .$

Further, the proof is similar to the proof of Theorem 4.1. \ \ \ \ \ \ \ \ \
\ \ \ \ \ \ \ \ \ \ \ \ \ \ \ \ \ \ \ \ \ \ \ \ \ \ \ \ \ \ \ \ $\square $%
\medskip

If $F$ satisfies the property $\alpha $ (resp. $\alpha ^{\prime }),$ we
obtain the following variant of Theorem 4.3.

\textbf{Theorem 4.4} \textit{Let }$X$ \textit{be a convex set}\ \textit{be
in a Hausdorff topological vector space }$E,$ $Y$\textit{\ be a compact set
in a Hausdorff topological vector space }$Z$\textit{\ and let }$S$\textit{\
be a pointed closed convex cone in }$Z$\textit{\ with its interior int}$%
S\neq \emptyset .$\textit{\ Let }$F:X\times X\rightrightarrows Y$\textit{\
be a set-valued map with non-empty values.}

\textit{i) Suppose that }$F$\textit{\ satisfies the property }$\alpha $. 
\textit{If }$F$ \textit{is transfer type-(iii) properly }$S-$\textit{%
quasi-concave in the first argument on }$X\times X$ \textit{and} \textit{also%
} \textit{satisfies the condition} $\gamma ,$\textit{\ then, there exist }$%
z_{1}\in $Max$\overline{\cup _{x\in X}F(x,x)}$\textit{\ and }$z_{2}\in $Min$%
\overline{\cup _{x\in X}\text{Max}_{w}F(x,X)}$ \textit{such that }$z_{1}\in
z_{2}+S.$

\textit{ii)} \textit{Suppose that that }$F$\textit{\ satisfies the property }%
$\alpha ^{\prime }$. \textit{If }$F$ \textit{is\ transfer type-(iii)
properly }$S-$\textit{quasi-convex\ in the second argument on }$X\times X$ 
\textit{and} \textit{also} \textit{satisfies the condition} $\gamma ^{\prime
},$\textit{\ then there exist }$z_{1}\in $Min$\overline{\cup _{x\in X}F(x,x)}
$\textit{\ and }$z_{2}\in $Max$\overline{\cup _{y\in X}\text{Min}_{w}F(X,y)%
\text{ }}$\textit{such that }$z_{1}\in z_{2}-S.$\medskip

We obtain the following corollary of Theorem 4.4, for the case of the
real-valued maps.

\textbf{Corollary 4.3} \textit{Let }$X$ \textit{be a convex set in a
Hausdorff topological vector space }$E,Y$\textit{\ be a compact set in }$%
\mathit{R}$\textit{\ and let }$S$\textit{\ be a pointed closed convex cone
in }$R$\textit{\ with its interior int}$S\neq \emptyset .$\textit{\ Let }$%
F:X\times X\rightrightarrows Y$\textit{\ be a set-valued map with non-empty
values. \newline
i) Suppose that }$\tbigcup_{y\in X}F(x,y)$\textit{\ is a compact set for
each }$x\in X$. \textit{If }$F$ \textit{is transfer type-(iii) properly }$S-$%
\textit{quasi-concave in the first argument on }$X\times X$ \textit{and} 
\textit{satisfies the condition} $\gamma ,$\textit{\ then, there exist }$%
z_{1}\in $Max$\overline{\cup _{x\in X}F(x,x)}$\textit{\ and }$z_{2}\in $Min$%
\overline{\cup _{x\in X}\text{Max}_{w}F(x,X)}$ \textit{such that }$z_{1}\in
z_{2}+S.$\textit{\ }

\textit{ii) Suppose that }$\tbigcup_{x\in X}F(x,y)$\textit{\ is a compact
set for each }$y\in X$. \textit{If }$F$ \textit{is transfer type-(iii)
properly }$S-$\textit{quasi-convex in the second argument on }$X\times X$ 
\textit{and} \textit{satisfies the condition} $\gamma ^{\prime },$\textit{\
then, there exist }$z_{1}\in $Min$\overline{\cup _{x\in X}F(x,x)}$\textit{\
and }$z_{2}\in $Max$\overline{\cup _{y\in X}\text{Min}_{w}F(X,y)}$ \textit{%
such that }$z_{1}\in z_{2}-S.$\medskip\ 

\textit{Example 4.3} Let $X=[0,1],$ $Y=[-1,1],$ $S=[0,\infty )$ and $%
F:X\times X\rightrightarrows Y$ be defined by $F(x,y)=\left\{ 
\begin{array}{c}
\lbrack 0,y]\text{ if }0\leq x\leq y\leq 1; \\ 
\lbrack -x,y]\text{ if }0\leq y<x\leq 1.%
\end{array}%
\right. $

We notice that $F$ is not continuous on $X.$

According to Examples 3.7 and 3.8, $F$ is transfer type-(iii)\textit{\ }%
properly $S-$quasi-concave in the first argument on $X\times X$ and it has
the property $\gamma .$

All the assumptions of Corollary 4.3 are fulfilled, then, there exists the
elements $z_{1}\in $Max$\overline{\cup _{x\in X}F(x,x)}$\ and $z_{2}\in $Min$%
\overline{\cup _{x\in X}\text{Max}_{w}F(x,X)}$ such that $z_{1}\in z_{2}+S.$

It is also true that:

$\cup _{x\in X}F(x,x)=\cup _{x\in X}[0,x]=[0,1];$ Max$\overline{\cup _{x\in
X}F(x,x)}=\{1\};$ Max$_{w}F(x,X)=\{1\}$ and Min$\overline{\cup _{x\in X}%
\text{Max}_{w}F(x,X)}=\{1\}.$

Then, taking $z_{1}=1$ and $z_{2}=1,$ we have $z_{1}\in z_{2}+S.$ \ \ \
\medskip

We introduce the following definition which concerns the convexity
properties of set-valued maps with two variables. It will be used to obtain
different minimax inequalities.

\textbf{Definition 4.1 }Let $X$ be a (n-1) dimensional simplex of a Hausdorff%
\textit{\ }topological vector space\textit{\ }$E,Y$ a subset of a Hausdorff
topological vector space\textit{\ }$Z$ and let $S$ be a pointed closed
convex cone in $Z$ with its interior int$S\neq \emptyset .$ Let $F:X\times
X\rightrightarrows Y$ be a set valued map with nonempty values.\newline
$F$ \textit{is weakly }$z-$\textit{convex\ }on $X$ for $z\in A\subseteq Z$,
iff for each $z\in A$ and $x_{1},...,x_{n}\in X,$ there exist $%
y_{1}^{z},y_{2}^{z},...,y_{n}^{z}\in X$ and $g^{z}\in C^{\ast }(\Delta
_{n-1})$ such that $F(x_{i},y_{i}^{z})\cap (z+S)\neq \emptyset $ for each $%
i\in \{1,2,...,n\}$ implies $F(\tsum_{i=1}^{n}\lambda
_{i}x_{i},\tsum_{i=1}^{n}g_{i}^{z}(\lambda _{i})y_{i}^{z})\cap (z+S)\neq
\emptyset $ for each $(\lambda _{1},\lambda _{2},...,\lambda _{n})\in \Delta
_{n-1}.$\medskip

\textit{Example} \textit{4.4} Let $X=[0,1],$ $Y=[0,1],$ $S=[0,\infty )$ and $%
F:X\times X\rightrightarrows Y$ be defined by $F(x,y)=\left\{ 
\begin{array}{c}
\lbrack 0,x]\text{ if }0\leq x\leq y\leq 1; \\ 
\lbrack 0,1]\text{ if }0\leq y<x\leq 1.%
\end{array}%
\right. $

For each $z\in \lbrack 0,1)$ and $x_{1},x_{2},...,x_{n}\in X,$ there exists $%
y_{1}^{z},y_{2}^{z},...,y_{n}^{z}\in X$ with $0\leq x_{i}\leq y_{i}^{z}$ for
each $i\in \{1,2,...,n\},$ such that $F(x_{i},y_{i}^{z})\cap
(z+S)=[0,x_{i}]\cap \lbrack z,\infty )\neq \emptyset $. It follows that $%
z\leq $min$_{i=1,...,n}\{x_{i}\}.$ Consequently, $z\leq
\tsum_{i=1}^{n}\lambda _{i}x_{i}$ and $0\leq x_{i}\leq
\tsum_{i=1}^{n}g_{i}^{z}(\lambda _{i})y_{i}^{z}$ for each $i\in
\{1,2,...,n\} $, $g^{z}\in C^{\ast }(\Delta _{n-1})$ and $\lambda =(\lambda
_{1},\lambda _{2},...,\lambda _{n})\in \Delta _{n-1}$. Then, $%
F(\tsum_{i=1}^{n}\lambda _{i}x_{i},\tsum_{i=1}^{n}g_{i}^{z}(\lambda
_{i})y_{i}^{z})=[0,\tsum_{i=1}^{n}\lambda _{i}x_{i}].$

Hence, $F(\tsum_{i=1}^{n}\lambda _{i}x_{i},\tsum_{i=1}^{n}g_{i}^{z}(\lambda
_{i})y_{i}^{z})\cap (z+S)=[0,\tsum_{i=1}^{n}\lambda _{i}x_{i}]\cap \lbrack
z,\infty )\neq \emptyset .$

For $z=1$ and for any $x_{1},x_{2},...,x_{n}\in X,$ there exists $%
y_{1}^{z},y_{2}^{z},...,y_{n}^{z}\in X$ with $0\leq y_{i}^{z}<x_{i}$ for
each $i\in \{1,2,...,n\},$ such that $F(x_{i},y_{i}^{z})\cap (z+S)=[0,1]\cap
\lbrack 1,\infty )\neq \emptyset $ for each $i\in \{1,2,...,n\}.$ We have
that $0\leq \tsum_{i=1}^{n}g_{i}^{z}(\lambda _{i})y_{i}^{z}<x_{i}$ for each $%
i\in \{1,2,...,n\}$, $g_{i}^{z}\in C^{\ast }(\Delta _{n-1})$ and $\lambda
=(\lambda _{1},\lambda _{2},...,\lambda _{n})\in \Delta _{n-1}$ and then, $%
F(\tsum_{i=1}^{n}\lambda _{i}x_{i},\tsum_{i=1}^{n}g_{i}^{z}(\lambda
_{i})y_{i}^{z})=[0,1].$ Therefore, $F(\tsum_{i=1}^{n}\lambda
_{i}x_{i},\tsum_{i=1}^{n}g_{i}^{z}(\lambda _{i})y_{i}^{z})\cap
(z+S)=[0,1]\cap \lbrack 1,\infty )\neq \emptyset .\medskip $

If $f$ is a mapping, we obtain the following definition.

\textbf{Definition 4.2 }Let $X$ be a (n-1) dimensional simplex of a\textit{\ 
}Hausdorff\textit{\ }topological vector space\textit{\ }$E,Y$ a subset of a
Hausdorff topological vector space\textit{\ }$Z$ and let $S$ be a pointed
closed convex cone in $Z$ with its interior int$S\neq \emptyset .$ Let $%
f:X\times X\rightarrow Y$ be a mapping.

$f$ \textit{is weakly }$z$-\textit{convex\ }on $X$ for $z\in A\subseteq Z$,
iff for each $z\in A$ and $x_{1},...,x_{n}\in X,$ there exist $%
y_{1}^{z},y_{2}^{z},...,y_{n}^{z}\in X$, $g^{z}\in C^{\ast }(\Delta _{n-1})$
such that $f(x_{i},y_{i}^{z})\in z+S$ for each $i\in \{1,2,...,n\}$ implies $%
f(\tsum_{i=1}^{n}\lambda _{i}x_{i},\tsum_{i=1}^{n}g_{i}^{z}(\lambda
_{i})y_{i}^{z})\in z+S$ for each $(\lambda _{1},\lambda _{2},...,\lambda
_{n})\in \Delta _{n-1}.$\medskip

\textit{Example} \textit{4.5} Let $X=[0,1],$ $Y=[0,1]\times \lbrack 0,1],$ $%
S=R_{+}^{2}$ and $f:X\times X\rightarrow Y$ be defined by $f(x,y)=\left\{ 
\begin{array}{c}
(x,y)\text{ if }0\leq x\leq y\leq 1; \\ 
(1,y)\text{ if }0\leq y<x\leq 1.%
\end{array}%
\right. $

For each $z=(z^{\prime },z^{\prime \prime })\in \lbrack 0,1)\times \lbrack
0,1]$ and $x_{1},x_{2},...,x_{n}\in X,$ there exists $%
y_{1}^{z},y_{2}^{z},...,$

$y_{n}^{z}\in X$ with $0\leq x_{i}\leq y_{i}^{z}$ for each $i\in
\{1,2,...,n\}$ such that $(x_{i},y_{i}^{z})=f(x_{i},y_{i}^{z})\in
(z+S)=[z^{\prime },\infty )\times \lbrack z^{\prime \prime },\infty ).$ It
follows that $z^{\prime }\leq $min$_{i=1,...,n}\{x_{i}\}.$ Consequently, $%
z^{\prime }\leq \tsum_{i=1}^{n}\lambda _{i}x_{i}$ and $0\leq x_{i}\leq
\tsum_{i=1}^{n}g_{i}^{z}(\lambda _{i})y_{i}^{z}$ for each $i\in
\{1,2,...,n\} $, $g^{z}\in C^{\ast }(\Delta _{n-1})$ and $\lambda =(\lambda
_{1},\lambda _{2},...,\lambda _{n})\in \Delta _{n-1}$ and then, $%
(\tsum_{i=1}^{n}\lambda _{i}x_{i},\tsum_{i=1}^{n}g_{i}^{z}(\lambda
_{i})y_{i}^{z})=f(\tsum_{i=1}^{n}\lambda
_{i}x_{i},\tsum_{i=1}^{n}g_{i}^{z}(\lambda _{i})y_{i}^{z})\in
(z+S)=[z^{\prime },\infty )\times \lbrack z^{\prime \prime },\infty ).$

For $z=(1,y)$ with $y\in \lbrack 0,1)$ and for any $x_{1},x_{2},...,x_{n}\in
X,$ there exists $y_{1}^{z},y_{2}^{z},...,y_{n}^{z}\in X$ with $0\leq
y_{i}^{z}<x_{i}$ for each $i\in \{1,2,...,n\}$ such that $%
(1,y_{i}^{z})=f(x_{i},y_{i}^{z})\in (z+S)=[1,\infty )\times \lbrack y,\infty
).$ We have that $0\leq \tsum_{i=1}^{n}g_{i}^{z}(\lambda
_{i})y_{i}^{z}<x_{i} $ for each $i\in \{1,2,...,n\}$, $g^{z}\in C^{\ast
}(\Delta _{n-1})$ and $\lambda =(\lambda _{1},\lambda _{2},...,\lambda
_{n})\in \Delta _{n-1}$ and then, $(1,\tsum_{i=1}^{n}g_{i}^{z}(\lambda
_{i})y_{i}^{z})=f(\tsum_{i=1}^{n}\lambda
_{i}x_{i},\tsum_{i=1}^{n}g_{i}^{z}(\lambda _{i})y_{i}^{z})\in
(z+S)=[1,\infty )\times \lbrack y,\infty ).\medskip $

Theorem 4.5 is a minimax theorem in which the set-valued map satisfies the
above defined property.

\textbf{Theorem} \textbf{4.5} \textit{Let }$X$\textit{\ be a (n-1)
dimensional simplex of a Hausdorff topological vector space }$E,$ $Y$\textit{%
\ a compact set in a Hausdorff topological vector space }$Z$\textit{\ and
let }$S$\textit{\ be a pointed closed convex cone in }$Z$\textit{\ with its
interior int}$S\neq \emptyset .$\textit{\ Let }$F:X\times X\rightrightarrows
Y$\textit{\ be a set-valued map with non-empty values such that }$\cup
_{x\in X}F(x,x)$ and $\cup _{y\in X}F(x,y)$\textit{\ are compact sets for
each }$x\in X$\textit{. Suppose the following conditions are fulfilled:%
\newline
(i) }$\mathit{F}$\textit{\ is weakly }$z-$\textit{convex\ for each }$z\in $%
Min$\overline{\cup _{x\in X}\text{Max}_{w}F(x,X)};$

\textit{(ii) for each }$x\in X,$\textit{\ }Min$\overline{\cup _{x\in X}\text{%
Max}_{w}F(x,X)}\subset F(x,X)-S.$

\textit{Then, }Min$\overline{\cup _{x\in X}\text{Max}_{w}F(x,X)}\subset $Max$%
\cup _{x\in X}F(x,x)-S.$\textit{\medskip }

\textit{Proof.} According to assumptions and Lemma 2.1, Max$_{w}F(x,X)\neq
\emptyset $ for each $x\in X$ and Min$\overline{\cup _{x\in X}\text{Max}%
_{w}F(x,X)}\neq \emptyset .$

Let $z\in $Min$\overline{\cup _{x\in X}\text{Max}_{w}F(x,X)}$ and let us
define the set-valued map $T:X\rightrightarrows X$ by $T(x)=\{y\in
X:F(x,y)\cap (z+S)\neq \emptyset \}$ for each $x\in X.$ According to
assumption (ii), it follows that $T(x)$ is nonempty for each $x\in X.$

According to Assumption (i), we have that $T$ is weakly naturally
quasi-convex: for any $x_{1},x_{2},...,x_{n}\in X$ and $z\in Y,$ there exist 
$y_{1}^{z},y_{2}^{z},...,y_{n}^{z}\in X$, $g^{z}\in C^{\ast }(\Delta
_{n-1}), $ such that, if $y_{i}\in T(x_{i})$ for each $i\in \{1,2,...,n\},$
then, $\tsum_{i=1}^{n}g_{i}^{z}(\lambda _{i})y_{i}^{z}\in
T(\tsum_{i=1}^{n}\lambda _{i}x_{i})$ for each $\lambda =(\lambda
_{1},\lambda _{2},...,\lambda _{n})\in \Delta _{n-1}.$

Therefore, according to the fixed point Theorem 2.1, there exists $x^{\ast
}\in T(x^{\ast }),$ that is, $F(x^{\ast },x^{\ast })\cap (z+S)\neq \emptyset
.$ Then, according to Lemma 2.1, we have $z\in F(x^{\ast },x^{\ast
})-S\subset \cup _{x\in X}F(x,x)-S\subset $Max$\cup _{x\in X}F(x,x)-S.$ \ \
\ \ \ \ \ \ \ \ \ \ \ \ \ \ \ \ \ \ \ \ \ \ \ \ \ \ \ \ \ \ \ \ \ \ \ \ \ \
\ \ \ \ \ \ \ \ \ \ $\square $\medskip

\textit{Example} \textit{4.6} Let $X=[0,1],$ $Y=[0,1],$ $S=[0,\infty )$ and $%
F:X\times X\rightrightarrows Y$ be defined by $F(x,y)=\left\{ 
\begin{array}{c}
\lbrack 0,x]\text{ if }0\leq x\leq y\leq 1; \\ 
\lbrack 0,1]\text{ if }0\leq y<x\leq 1.%
\end{array}%
\right. $

We saw in Example 4.2 that $F$ is weakly $z-$convex for each $z\in Z.$

Further, we have that, $\cup _{x\in X}F(x,x)=\cup _{x\in X}[0,x]=[0,1]$ and
for each $x\in X,$ $\cup _{y\in X}F(x,y)=[0,1]$, so that, $\cup _{x\in
X}F(x,x)$ and $\cup _{y\in X}F(x,y)$ are compact sets, for each $x\in X.$

It is also true that Max$_{w}F(x,X)=\{1\}$ and Min$\overline{\cup _{x\in X}%
\text{Max}_{w}F(x,X)}=\{1\}.$

$F(x,X)-S=(-\infty ,1]$ and then, for each\textit{\ }$x\in X,$\textit{\ }Min$%
\overline{\cup _{x\in X}\text{Max}_{w}F(x,X)}$

$\subset F(x,X)-S.$ All the assumptions of Theorem 4.4 are fulffiled.

Then, $\{1\}=$Min$\overline{\cup _{x\in X}\text{Max}_{w}F(x,X)}\subset $Max$%
\cup _{x\in X}F(x,x)-S=(-\infty ,1].\medskip $

The next corollary is obtained by considering single valued mappings, as a
particular case, in Theorem 4.5.

\textbf{Corollary 4.4} \textit{Let }$X$\textit{\ be an (n-1) dimensional
simplex of a Hausdorff topological vector space }$E,Y$\textit{\ a compact
set in a Hausdorff topological vector space }$Z$\textit{\ and let }$S$%
\textit{\ be a pointed closed convex cone in }$Z$\textit{\ with its interior
int}$S\neq \emptyset .$\textit{\ Let }$f:X\times X\rightarrow Y$\textit{\ be
a mapping such that }$\cup _{y\in X}f(x,y)$\textit{\ and }$\cup _{x\in
X}f(x,x)$\textit{\ are compact sets for each }$x\in X$\textit{. Suppose the
following conditions are fulfilled:}

\textit{(i) }$\mathit{f}$\textit{\ is weakly }$z-$\textit{convex\ for each }$%
z\in $Min$\overline{\cup _{x\in X}\text{Max}_{w}f(x,X)};$

\textit{(ii) for each }$x\in X,$\textit{\ }Min$\overline{\cup _{x\in X}\text{%
Max}_{w}f(x,X)}\subset f(x,X)-S.$

\textit{Then, }Min$\overline{\cup _{x\in X}\text{Max}_{w}f(x,X)}\subset $Max$%
\cup _{x\in X}f(x,x)-S.$\textit{\medskip }

\textit{Example} \textit{4.7} Let $X=[0,1],$ $Y=[0,1]\times \lbrack 0,1],$ $%
S=IR_{+}^{2}$ and $f:X\times X\rightarrow Y$ be defined by $f(x,y)=\left\{ 
\begin{array}{c}
(x,y)\text{ if }0\leq x\leq y\leq 1; \\ 
(1,1)\text{ if }0\leq y<x\leq 1.%
\end{array}%
\right. $

We notice that $f$ is not continuous. The mapping $f$ is weakly $z-$convex
for each $z\in Y$. According to the definition of $f$, $\cup _{y\in
X}f(x,y)=\{x\}\times \lbrack x,1]\cup \{(1,1)\}$ and $\cup _{x\in
X}f(x,x)=\{(x,x):x\in \lbrack 0,1]\},$ which are compact sets.

The following equalities take place:

Max$_{w}\cup _{y\in X}f(x,y)=\{1\}\times \lbrack 0,1]\cup \lbrack 0,1]\times
\{1\}$ and

Min$\overline{\cup _{x\in X}\text{Max}_{w}f(x,X)}=\{1\}\times \lbrack
0,1]\cup \lbrack 0,1]\times \{1\}.$

Finally, we have that, for each\textit{\ }$x\in X,$\textit{\ }Min$\overline{%
\cup _{x\in X}\text{Max}_{w}f(x,X)}\subset f(x,X)-S$ and then, all the
assumptions of the Corrollary are satisfied.

Hence, Min$\overline{\cup _{x\in X}\text{Max}_{w}f(x,X)}\subset $Max$\cup
_{x\in X}f(x,x)-S.$\textit{\medskip }

Another result is obtained in the same context of Theorem 4.5.

\textbf{Theorem} \textbf{4.6} \textit{Let }$X$\textit{\ be an (n-1)
dimensional simplex of a Hausdorff topological vector space }$E,Y$\textit{\
a compact set in a Hausdorff topological vector space }$Z$\textit{\ and let }%
$S$\textit{\ be a pointed closed convex cone in }$Z$\textit{\ with its
interior int}$S\neq \emptyset .$\textit{\ Let }$F:X\times X\rightrightarrows
Y$\textit{\ be a set-valued map with non-empty values such that }$\cup
_{x\in X}F(x,x)$ and $\cup _{y\in X}F(x,y)$\textit{\ are compact sets for
each }$x\in X$\textit{. Suppose the following conditions are fulfilled:%
\newline
(i) }$\mathit{F}$\textit{\ is weakly }$z-$\textit{convex\ for each }$z\in $%
Max$\overline{\cup _{y\in X}\text{Min}_{w}F(X,y)};$

\textit{(ii) for each }$x\in X,$\textit{\ }Max$\overline{\cup _{y\in X}\text{%
Min}_{w}F(X,y)}\subset F(X,y)+S.$

\textit{Then, }Max$\overline{\cup _{y\in X}\text{Min}_{w}F(X,y)}\subset $Min$%
\cup _{x\in X}F(x,x)+S.$\textit{\medskip }

\textit{Proof.} According to the assumptions and Lemma 2.1, Min$%
_{w}F(X,y)\neq \emptyset $ for each $y\in X$ and Max$\overline{\cup _{y\in X}%
\text{Min}_{w}F(X,y)}.$

Let $z\in $Max$\overline{\cup _{y\in X}\text{Min}_{w}F(X,y)}$ and let us
define the set-valued map $Q:X\rightrightarrows X$ by $Q(y)=\{x\in
X:F(x,y)\cap (z-S)\neq \emptyset \}$ for each $y\in X.$ According to the
assumption (ii), it follows that $Q(y)$ is non-empty for each $y\in X.$

Accordint to the Assumption (i), we have that $Q$ is weakly naturally
quasi-convex: for any $y_{1},y_{2},...,y_{n}\in X$ and $z\in Y,$ there exist 
$x_{1}^{z},x_{2}^{z},...,x_{n}^{z}\in X$ and $g^{z}\in C^{\ast }(\Delta
_{n-1}),$ such that, if $x_{i}\in Q(y_{i})$ for each $i\in \{1,2,...,n\},$
then, $\tsum_{i=1}^{n}g_{i}^{z}(\lambda _{i})x_{i}^{z}\in
Q(\tsum_{i=1}^{n}\lambda _{i}y_{i})$ for each $\lambda =(\lambda
_{1},\lambda _{2},...,\lambda _{n})\in \Delta _{n-1}.$

Therefore, according to the fixed point Theorem 2.1, there exists $y^{\ast
}\in Q(y^{\ast }),$ that is, $F(y^{\ast },y^{\ast })\cap (z-S)\neq \emptyset
.$ According to Lemma 2.1, we have that $z\in F(y^{\ast },y^{\ast
})+S\subset \cup _{x\in X}F(x,x)+S\subset $Min$\cup _{x\in X}F(x,x)+S.$ \ \
\ \ \ \ \ \ \ \ \ \ \ \ \ \ \ \ \ \ \ \ \ \ \ \ \ \ \ \ \ \ \ \ \ \ \ \ \ \
\ \ \ \ \ \ \ \ \ \ \ $\square $\medskip

The last result from this paper is stated now.

\textbf{Corollary} 4.5 \textit{Let }$X$\textit{\ be an (n-1) dimensional
simplex of a Hausdorff topological vector space }$E,Y$\textit{\ a compact
set in a Hausdorff topological vector space }$Z$\textit{\ and let }$S$%
\textit{\ be a pointed closed convex cone in }$Z$\textit{\ with its interior
int}$S\neq \emptyset .$\textit{\ Let }$f:X\times X\rightarrow Y$\textit{\ be
a mapping such that }$\cup _{y\in X}f(x,y)$\textit{\ and }$\cup _{x\in
X}f(x,x)$\textit{\ are compact sets for each }$x\in X$\textit{. Let us
suppose that the following conditions are fulfilled:}

\textit{(i) }$\mathit{f}$\textit{\ is weakly }$z-$\textit{convex\ for each }$%
z\in $Max$\overline{\cup _{y\in X}\text{Min}_{w}f(X,y)};$

\textit{(ii) for each }$x\in X,$\textit{\ }Max$\overline{\cup _{y\in X}\text{%
Min}_{w}f(X,y)}\subset f(X,y)+S.$

\textit{Then, }Max$\overline{\cup _{y\in X}\text{Min}_{w}f(X,y)}\subset $Min$%
\cup _{x\in X}f(x,x)+S.$\textit{\medskip }

\textbf{Concluding Remarks}

We have proven the existence of equilibria in minimax inequalities without
assuming any form of continuity of functions or set-valued maps. New
conditions of convexity have been introduced. The main tools to prove our
results have been a fixed-point theorem for weakly naturally quasi-concave
set valued maps and a constant selection for quasi-convex set-valued maps.
Several examples have been provided in order to illustrate our results.%
\textit{\medskip }

\textbf{REFERENCES}\textit{\medskip \medskip }

[1] Chen G. Y. (1991) A generalized section theorem and a minimax inequality
for a vector-valued mapping. Optimization \textbf{22}, 745-754

[2] Chuang C.-S., Lin L.-J. (2012): New existence theorems for
quasi-equilibrium problems and a minimax theorem on complete metric spaces.
J. Glob. Optim., DOI 10.1007/s10898-012-0004-3

[3] Ding X. and Yiran He (1998) Best Approximation Theorem for Set-valued
Mappings without Convex Values and Continuity. Appl Math. and Mech. English
Edition \textbf{19}(9), 831-836

[4] Fan K. (1953) Minimax theorems. Proceedings of the National Academy of
Sciences of the USA, 39, 42-47

[5] Fan K. (1972) A Minimax Inequality and Applications. Inequalities III,
Academic Press, New York, NY

[6] Fan K. (1964): Sur un th\'{e}or\`{e}me minimax. C. R. Acad. Sci. Paris
259, 3925-3928

Ferro F. (1989) A minimax theorem for Vector-Valued functions. J. Optim.
Theory Appl. 60, 19-31

[7] Ferro F. (191) A Minimax Theorem for Vector-Valued Functions. J. Optim.
Theory Appl. 68, 35--48

[8] Georgiev P. G. and TanakaT. (2000) Ky Fan's Inequality for Set-Valued
Maps with Vector-Valued Images. Nonlinear analysis and convex analysis
(Kyoto, 2000) S\={u}rikaisekikenky\={u}sho K\={o}ky\={u}roku No. 1187,
143-154 (2001)

[9] Ha C. W. (1980) Minimax and fixed point theorems. Mathematische Annalen
248, 73-77

[10] Jahn J. (2004) Vector Optimization: Theory, Applications and
Extensions. Springer, Berlin

[11] Kuroiowa D., Tanaka T. and Ha T. X. D. (1997) On cone convexity of
set-valued maps. Non-linear Analysis: Theory, Methods and Applications, 
\textbf{30}(3)1487-1496

[12] Li S. J., Chen G. Y. and Lee G. M. (2000) Minimax theorems for
set-valued mappings. J. Optim. Theory Appl. 106, 183-199

[13] Li S. J., Chen G. Y., Teo K. L. and Yang X. Q. (2003) Generalized
minimax inequalities for set-valued mappings. J. Math. Anal. Appl. 281,
707-723

[14] Li Z. F. and \ Wang S. Y. (1998) A type of minimax inequality for
vector-valued mappings. J Math Anal Appl. 227, 68--80

[15] Lin Y. J., Tian G. (1993) Minimax inequalities equivalent to the
Fan-Knaster-Kuratowski-Mazurkiewicz theorems. Appl. Math. Opt. \textbf{28}%
(2), 173-179

[16] Lin L.-J., Du W.-S. (2008) Systems of equilibrium problems with
applications to new variants of Ekeland's variational principle, fixed point
theorems and parametric optimization problems. J. Glob. Optim. \textbf{40}%
(4), 663-677

[17] D.T. Luc D.T., Vargas C. A. (1992) A seddle point theorem for
set-valued maps. Nonlinear Anal. 18, 1-7

[18] Luo X. Q. (2009) On some generalized Ky Fan minimax inequalities. Fixed
Point Theory Appl 2009, 9 . Article ID 194671doi:10.1155/2009/194671

[19] Nessah R. and G. Tian G. (2013) Existence of Solution of Minimax
Inequalities, Equilibria in Games and Fixed PointsWithout Convexity and
Compactness Assumptions. J Optim Theory Appl, 157, 75--95 DOI
10.1007/s10957-012-0176-5

[20] Nieuwenhuis J. W. (1983) Some minimax theorems in vector-valued
functions. J. Optim. Theory Appl. 40, 463-475

[21] Patriche M. Fixed point theorems and applications in theory of games.
Fixed Point Theory, forthcoming

[22] Sion M. (1958) On General Minimax Theorems. Pac. J. of Math. 8, 171--176

[23] Smajdor W. (2000) Set-valued quasiconvex functions and their constant
selections. Functional Equations and Inequalities, Kluwer Academic
Publishers (2000)

[24] Tanaka T. (1998) Some Minimax Problems of Vector-Valued Functions. J.
Optim. Theory Appl. 59, 505--524

[25] Tanaka T. (1994) Generalied quasiconvexities, cone saddle points and
minimax theorems for vector-valued functions. J. Optim. Theory Appl. 81,
355-377

[27] Tuy H. (2011) A new topological minimax theorem with application. J.
Glob. Optim. \textbf{50(}3), 371-378

[28] Zhang Q. B. and Cheng C. Z. (2007) Some fixed-point theorems and
minimax inequalities in FC-space. J. Math. Anal. Appl. 328, 1369-1377

[29] Zhang Q. B., Liu M. J. and Cheng C. Z. (2009) Generalized saddle points
theorems for set-valued mappings in locally generalized convex spaces.
Nonlinear Anal., 71, 212-218

[30] Zhang Q., Cheng C. and Li X. (2013) Generalized minimax theorems for
two set-valued mappings. J. Ind. Man. Optim. 9 1, 1-12,
doi:10.3934/jimo.2013.9.1

[31] Zhang Y., Li S. J. (2012) Minimax theorems for scalar set-valued
mappings with nonconvex domains and applications. J. Glob. Optim. (2012) DOI
10.1007/s10898-012-9992-2

[32] Zhang Y. and Li S. J.: Ky Fan minimax inequalities for set-valued
mappings. Fixed Point Theory Appl. 2012, 2012:64

\end{document}